%% file: ModelComparison.tex
\documentclass[11pt]{amsart}
\input{PreambleSep2020.tex}

\title[$2$-complicial sets and $\Theta_2$-spaces]{An explicit comparison between\\ $2$-complicial sets and $\Theta_2$-spaces}

\author{Julia E. Bergner}
\address{Department of Mathematics, University of Virginia, Charlottesville,
USA}
\email{jeb2md@virginia.edu}

\author{Viktoriya Ozornova}
\address{Max Planck Institute for Mathematics, Bonn, Germany}
\email{viktoriya.ozornova@mpim-bonn.mpg.de}

\author{Martina Rovelli}
\address{Department of Mathematics and Statistics, 
University of Massachusetts, 
Amherst,
USA
}
\email{rovelli@math.umass.edu} 

\makeatletter
\@namedef{subjclassname@2020}{\textup{2020} Mathematics Subject Classification}
\makeatother

\begin{document}
\begin{abstract}
We produce a direct Quillen equivalence between two models of $(\infty,2)$-categories: the complete Segal $\Theta_2$-spaces due to Rezk and the $2$-complicial sets due to Verity.
\end{abstract}

\maketitle

\section*{Introduction}

The language of higher categories provides a way to describe many phenomena in areas of mathematics as diverse as topology, algebra, geometry, and mathematical physics.  In a higher categorical structure, we not only have functions between objects, but functions between those functions and possibly further iterations of this idea, encoded by the notion of a $k$-morphism between $(k-1)$-morphisms.  One might initially assume that these higher morphisms should satify conditions like associativity in the usual way, but for many natural examples they only hold up to isomorphism or, in topological settings, up to homotopy.  In the latter situation, it is convenient to work in the setting of $(\infty, n)$-categories, in which we have $k$-morphisms for arbitrarily large $k$, but they are all weakly invertible for $k>n$.  These higher invertible morphisms provide a means for conveniently encoding the ``up to isomorphism'' data in the lower morphisms.

There have been many different approaches to realizing $(\infty,n)$-categories as concrete mathematical objects; such realizations are often called \emph{models} for $(\infty,n)$-categories.  A natural question, then, is whether these different models really do encode the same information, namely, whether we can establish an appropriate equivalence between them.  Much work has been done in this direction, but there are still proposed models for which we do not have such comparisons.  In some other cases, we know by general results that models must be equivalent, but do not have an explicit equivalence. 

The motivation for this paper is to give an explicit comparison between two of these models, the complete Segal $\Theta_n$-spaces as defined by Rezk \cite{rezkTheta} and the $n$-complicial sets as defined by Verity \cite{VerityComplicialI,EmilyNotes,or}; we give such a comparison when $n=2$, for which more tools are available. Let us give a brief description of these two models.

A complete Segal $\bT$-space is described by a diagram of spaces indexed by 2-categories freely generated by pasting diagrams such as
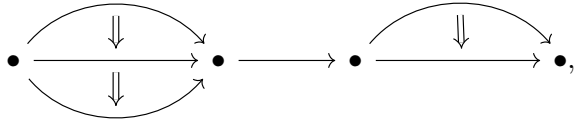
\[
\begin{tikzcd}[row sep=3.2cm, column sep=2.2cm]
  \bullet \arrow[r, bend left=50, "", ""{name=U,inner sep=2pt,below}]
  \arrow[r, ""{near end, xshift=0.2cm}, ""{name=D,inner sep=2pt},""{name=M,inner sep=2pt, below}]
  \arrow[r, bend right=50, ""{below}, ""{name=DD,inner sep=2pt}]
  & \bullet\arrow[r, ""]
    & [-0.9cm]
    %& 
    \bullet\arrow[r, bend left=50, "", ""{name=U1,inner sep=2pt,below}]
  \arrow[r, ""{near end, xshift=0.2cm}, ""{name=D1,inner sep=2pt, xshift=0.05cm},""{name=M1,inner sep=2pt, below}]
%  \arrow[r, bend right=50, "h"{below}, ""{name=DD1,inner sep=2pt}]
  & \bullet,
    \arrow[Rightarrow, from=U, to=D, ]
    \arrow[Rightarrow, from=M, to=DD, ""]  %{near start}
  \arrow[Rightarrow, from=U1, to=D1, ""]
\end{tikzcd}\]
which the expert reader may recognize as the generic element of Joyal's \emph{cell category} $\bT$.
In contrast, a $2$-complicial set is given by a simplicial set (with a suitable marking) in which a $k$-simplex represents a diagram indexed by a truncated \emph{oriental}, which is a free $2$-category generated by a standard simplex, such as
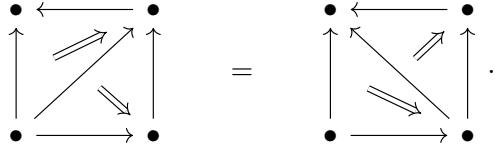
\[
\simpfverA{}{}{}{}{}{}{}{}{}\simpfverAcontinued{}{\bullet}{\bullet}{\bullet}{\bullet}.
\] 

A common way to show that two models are equivalent is to show that appropriate model categories for each are Quillen equivalent to each other.  In this paper, we seek to establish such a Quillen equivalence between the model structure $\spsh{\bT}_{p,(\infty,2)}$ for complete Segal $\bT$-spaces and the model structure $\msset_{(\infty,2)}$ for 2-complicial sets.

Combining several prior results by different groups of authors, we already know that the two model categories are Quillen equivalent via a rather lengthy zig-zag of Quillen equivalences between different models. 
Although we do not expect the reader to be familiar with all these models of $(\infty,2)$-categories, to give an idea of the complexity of the comparison here is a diagram of an essentially optimal zig-zag of Quillen equivalences, extracted from \cite{GHL}:
\begin{equation}
\label{DiagramIntro}
    \begin{tikzpicture}[baseline=-2.2cm]
    \def\r{3.7cm}
\node[outer sep=0.2cm] (A1) at (0,-0.4*\r) {$\spsh{\bT}_{i,(\infty,2)}$};
\node[outer sep=0.2cm] (A2) at (0,0) {$\spsh{\bT}_{p,(\infty,2)}$};
\node (A3) at (0,-0.8*\r) {$\spsh{(\Delta\times\Delta)}_{(\infty,2)}$};
\node (A4) at (0,-1.2*\r) {$P\cat(\spsh{\Delta})_{(\infty,2)}$};
\node (A5) at (0.9*\r,-1.2*\r) {$\vcat{\spsh{\Delta}_{(\infty,1)}}$};
\node (B1) at (1.8*\r, -1.2*\r) {$\vcat{\psh{\Delta}_{(\infty,1)}}$.}; 
\node (C1) at (1.8*\r, -0.8*\r) {$\vcat{\sset^+_{(\infty,1)}}$}; 
\node (C2) at (1.8*\r, -0.4*\r) {$\sset^{sc}_{(\infty,2)}$};
\node (D1) at (1.8*\r,0) {$\msset_{(\infty,2)}.$};

 \path (A1) -- node(comp1)[sloped,transform shape,inner sep=0pt,xscale=2.5] {$\simeq$} 
        node[font=\tiny, right of=comp1]{\cite{rezk}}(A2);
 \path (A1) -- node(comp2)[sloped,transform shape,inner sep=0pt,xscale=2.5, rotate=180] {$\simeq$} 
        node[font=\tiny, right of=comp2]{\cite{br2}}(A3);
\path (A3) -- node(comp3)[sloped,transform shape,inner sep=0pt,xscale=2.5, rotate=180] {$\simeq$} 
        node[font=\tiny, right of=comp3]{\cite{br2}}(A4);   
%%%%
\path (A4) -- node(comp4)[sloped,transform shape,inner sep=0pt,xscale=2.5, outer sep=0pt] {$\simeq$} 
        node[font=\tiny, below of=comp4, yshift=0.4cm]{\cite{br1}}(A5);   
\path (A5) -- node(comp5)[sloped,transform shape,inner sep=0pt,xscale=2.5, outer sep=0pt] {$\simeq$} 
        node[font=\tiny, below of=comp4, yshift=0.4cm]{\cite{JT}}(B1); 
\path (B1) -- node(comp6)[sloped,transform shape,inner sep=0pt,xscale=2.5, rotate=180] {$\simeq$} 
        node[font=\tiny, left of=comp6]{\cite{htt}}(C1);
\path (C1) -- node(comp7)[sloped,transform shape,inner sep=0pt,xscale=2.5, rotate=180] {$\simeq$} 
        node[font=\tiny, left of=comp7]{\cite{lurieGoodwillie}}(C2);
\path (C2) -- node(comp8)[sloped,transform shape,inner sep=0pt,xscale=2.5, rotate=180] {$\simeq$} 
        node[font=\tiny, left of=comp7]{\cite{GHL}}(D1);
    \end{tikzpicture}
\end{equation}

To simplify the comparison, the goal of this paper is to produce the following direct Quillen equivalence.

\begin{unnumberedtheorem} 
There is a Quillen equivalence between complete Segal $\bT$-spaces, presented by the model category $\spsh{\bT}_{p,(\infty,2)}$, and $2$-complicial sets, presented by the model category $\msset_{(\infty,2)}$.
\end{unnumberedtheorem}

In addition to providing a more transparent comparison between the two models, this direct comparison facilitates the transport of constructions between these model structures.  For example, the structure of $\bT$ makes the description of duals straightforward in $\bT$-spaces, whereas the join construction has been described for $2$-complicial sets by Verity \cite{VerityComplicialI}; this direct equivalence enables us to understand the analogous constructions in the less-convenient context for each construction.  

Let us now describe the main ingredients of the proof of our main theorem.
\begin{enumerate}[label=(\roman*),leftmargin=*, widest=ii]
    \item We use the compatibility of the $2$-categorical nerve valued in marked simplicial sets established by second and third authors in \cite{ORfundamentalpushouts}, to construct a left Quillen functor
    \[L\colon\spsh{\bT}_{(\infty,2)}\to\msset_{(\infty,2)}.\]
     
    \item To show that this left Quillen functor is in fact a Quillen equivalence, we use a result of Barwick and Schommer-Pries \cite{BarwickSchommerPries} to reduce the problem to showing that it preserves cells in dimensions $0$, $1$, and $2$.  In \cref{Sec:RecognizingCells} we use the intermediate comparisons of models from the diagram above, to identify these cells in each model and thereby show that $L$ does indeed preserve cells.
\end{enumerate}

The outline of the paper is as follows.  In \cref{GlobularSimplicial} we recall some necessary results about model categories of $2$-categories, $\bT$-spaces, and simplicial sets with marking, and of functors between them, such as suspension and nerve functors. In \cref{ComparisonStructure} we construct the adjunction between $\bT$-spaces and simplicial sets with marking that we and show that it is a Quillen pair.  We then describe how it follows from \cite{BarwickSchommerPries} that this adjunction is indeed a Quillen equivalence, modulo an explicit identification of the cells in the two models. In \cref{Sec:RecognizingCells} we then provide the desired identification of the cells in the two models.

\addtocontents{toc}{\protect\setcounter{tocdepth}{1}}
  \subsection*{Acknowledgements}
  The authors would like to thank Lennart Meier for helpful conversations on this project. The first-named author was partially supported by NSF grant DMS-1906281. The second-named author thankfully acknowledges the financial support by the DFG grant OZ 91/2-1 with the project nr.~442418934. This material is based upon work supported by the National Science Foundation under Grant No.\ DMS-1440140 while the second- and third-named authors were in residence at the Mathematical Sciences Research Institute in Berkeley, California, during the Spring 2020 semester.
  
  \newpage
\tableofcontents

\addtocontents{toc}{\protect\setcounter{tocdepth}{2}}
\section{Models of $(\infty,2)$-categories}
\label{GlobularSimplicial}

We assume the reader to be familiar with the basics of strict $2$-category theory (see e.g.~\cite{BorceuxHandbook1}) and with the language of model categories (see e.g.~\cite{Hirschhorn,hovey}), and we now recall some further preliminary material that we need in this paper. 

\subsection{Strict $2$-categories}

The category $2\cat$ of $2$-categories is defined as the category whose objects are (small) categories enriched over the category $\cat$ of $1$-categories.  In particular, a \emph{$2$-category} $\cD$ consists of a set of objects and for any objects $x,x'$ a $1$-category $\Map_{\cD}(x,x')$, together with a horizontal composition that defines a functor of hom-categories
$\circ\colon\Map_{\cD}(x,x')\times\Map_{\cD}(x',x'')\to\Map_{\cD}(x,x'').$

We consider the following model structure on $2\cat$ that was constructed by Lack in \cite[Thm~3.3]{lack1} (with a correction in \cite[Thm~4]{lack2}).

\begin{thm}
The category $2\cat$ of $2$-categories supports a model structure in which
\begin{itemize}[leftmargin=*]
    \item all $2$-categories are fibrant, and
    \item the weak equivalences are precisely the biequivalences of $2$-categories.
\end{itemize}
 \end{thm}

An important source of examples of $2$-categories is given by suspending $1$-categories, as follows.

\begin{defn} \label{suspension}
Let $\cD$ be a $1$-category $\cD$. The \emph{suspension} of $\cD$ is the $2$-category $\Sigma\cD$ in which
\begin{enumerate}[label=(\alph*), leftmargin=*, widest=c]
    \item there are two objects $x_{\bot}$ and $x_{\top}$; 
    \item the hom-$1$-categories given by 
\[
\Map_{\Sigma\cD}(a,b):=
\left\{
\begin{array}{cll}
\cD&\mbox{ if }a=x_{\bot},b=x_{\top}\\
{[0]}& \mbox{ if }a=b, \\
\varnothing& \mbox{ if } a=x_{\top},b=x_{\bot}; 
\end{array}
\right.
\]
and
\item there is no nontrivial horizontal composition.
\end{enumerate} 
This construction extends to a functor $\Sigma\colon\cat\to2\cat_{*,*}$ valued in the category of bipointed categories, namely categories endowed with a pair of (possibly equal) specified objects, and basepoint-preserving functors.
\end{defn}

The $2$-categorical suspension $\Sigma\cD$ appears in \cite{BarwickSchommerPries} as $\sigma(\cD)$.   It also often appears in the literature as a special case of a simplicial suspension. For instance, applying the nerve to hom categories of the suspension $\Sigma\cD$ gives a simplicial category $N_*(\Sigma\cD)$ that agrees with what was denoted by $U(N\cD)$ in \cite{bergner}, as $S(N\cD)$ in \cite{Joyal2007}, as $[1]_{N\cD}$ in \cite{htt}, and as ${\mathbbm{2}}[N\cD]$ in \cite{RiehlVerityNcoh}. 

\begin{notn}
We record the notation for the following (non-disjoint) families of $2$-categories.
\begin{itemize}[leftmargin=*]
    \item For $m\ge-1$, we denote by $[m]$ the finite ordinal with $m+1$ elements.
    
    \item For $j=0,1,2$, we denote by $C_j$ the free $j$-cell. These $2$-categories can be pictured as
    
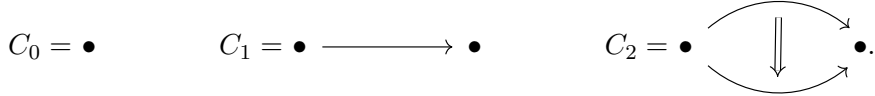
\[
    \begin{tikzcd}[row sep=3.2cm, column sep=1.8cm, inner sep=2pt, outer sep=1pt]
    C_0=&[-2cm] \bullet 
  &[-0.5cm] C_1= &[-2cm]\bullet\arrow[r]
    & %[-0.9cm]
    \bullet
&[-0.5cm] C_2=  &[-2cm] \bullet\arrow[r, bend left=40, "", ""{name=U,inner sep=2pt,below}, end anchor={[xshift=0.1cm]}]
  \arrow[r, bend right=40, "", ""{name=D,inner sep=2pt,above}, end anchor={[xshift=0.1cm]}]
  & \bullet.
    \arrow[Rightarrow, from=U, to=D]
     \end{tikzcd}
    \]

    \item For $m\ge0$ and $k_1,\ldots, k_m\geq 0$, we denote by $[m|k_1,\dots,k_m]$ the object of Joyal's cell category $\bT$, namely the full subcategory $\bT$ of $2\cat$ from \cite{JoyalDisks}.
        \item We denote by $\mathbb I$ the free-living isomorphism category. This category can be pictured as
        
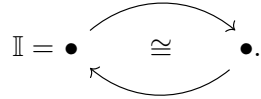
\[
            \begin{tikzcd}[row sep=3.2cm, column sep=1.8cm, inner sep=2pt, outer sep=1pt]
 \mathbb{I}=  &[-2cm] \bullet\arrow[r, bend left=40, "", ""{name=U,inner sep=2pt,below}, end anchor={[xshift=0.1cm]}, start anchor={[xshift=-0.05cm]}]
  & \bullet.\arrow[l, bend left=40, "", ""{name=D,inner sep=2pt,above}, end anchor={[xshift=-0.05cm]}, start anchor={[xshift=0.05cm]}]
    \arrow[phantom, "\cong", from=U, to=D]
     \end{tikzcd}
        \]
\end{itemize}
\end{notn}

%%%

\subsection{Complete Segal spaces as a model for $(\infty,1)$-categories}
\label{segalspacesection} We briefly recall the theory of complete Segal spaces, as first defined by Rezk in \cite{rezk}, of which the next model we discuss for $(\infty,2)$-categories is a generalization.

First, consider functors $X \colon \Deltaop \rightarrow \sset$.  For any $n \geq 1$, consider the \emph{Segal map}
\[ X_n \rightarrow \underbrace{X_1 \underset{X_0}{\times} \cdots \underset{X_0}{\times} X_1}_n \]
induced by the inclusion
\[ \underbrace{\Delta[1] \aamalg{\Delta[0]} \Delta[1]\aamalg{\Delta[0]}\ldots \aamalg{\Delta[0]} \Delta[1] }_{n}\rightarrow \Delta[n]\]
of the spine of the $n$-simplex into the $n$-simplex $\Delta[n]$.

\begin{defn}
A \emph{Segal space} is an injectively fibrant functor $X \colon \Deltaop \rightarrow \sset$ such that the Segal maps are weak equivalences of simplicial sets for all $n \geq 1$.
\end{defn}

The idea is that a Segal space behaves something like a category, with simplicial sets of objects and morphisms, but with composition defined only up to homotopy.  

However, to have a model for $(\infty, 1)$-categories, we do not want a simplicial set of objects (as in an internal category), but instead a discrete set of objects.  The most straightforward way to get such a model is to ask for the simplicial set $X_0$ to be discrete.

\begin{defn}
A \emph{Segal precategory} is a functor $X \colon \Deltaop \rightarrow \sset$ such that $X_0$ is a discrete simplicial set. We denote by $P\cat$ the full subcategory of $\spsh{\Delta}$ spanned by all Segal precategories.
A \emph{Segal category} is a Segal precategory that is also a Segal space.
\end{defn}

\begin{thm}[{\cite[Thm 5.1]{bergner3models}, \cite[Thm 6.4.4]{pelissier}}]
\label{secatmc}
The category $P\cat$ of Segal precategories admits a model structure in which
\begin{itemize}[leftmargin=*]
    \item the fibrant objects are the Segal categories; and
    \item the cofibrations are the monomorphisms.
\end{itemize}
We denote this model structure by $P\cat_{(\infty,1)}$.
\end{thm}

However, from the point of view of homotopy theory, asking for discreteness is awkward.  The completeness condition that we now describe can be more convenient from this perspective.

Let $N\mathbb I$ denote the nerve of the groupoid $\mathbb I$, and denote by $X_{\heq}$ the simplicial set $\Map(N\mathbb I, X)$, which is sometimes called the space of \emph{homotopy equivalences} of $X$. The unique map $N\mathbb I \rightarrow \Delta[0]$ induces a map 
\[ X_0 \rightarrow X_{\heq}. \]

\begin{defn}
A Segal space is \emph{complete} if this map $X_0 \rightarrow X_{\heq}$ is a weak equivalence of simplicial sets.
\end{defn}

Rezk builds a supporting model structure for the homotopy theory of complete Segal spaces.
 
\begin{thm}[{\cite[Thm 7.2]{rezk}}]
\label{Rezkmodelstructure}
The category $\sset^{\Deltaop}$ of simplicial spaces admits a model structure in which
\begin{itemize}[leftmargin=*]
    \item the fibrant objects are the complete Segal spaces, and
    
    \item the cofibrations are the monomorphisms.
\end{itemize}
We denote this model structure by $\sset^{\Deltaop}_{(\infty,1)}$.
\end{thm}

This model structure can be obtained
by taking the left Bousfield localization of the injective model structure on $\spsh{\Delta}$ with respect to the following set of maps:
\begin{enumerate}[leftmargin=*]
    \item the \emph{Segal acyclic cofibrations}
    \[ \underbrace{\Delta[1] \aamalg{\Delta[0]} \Delta[1]\aamalg{\Delta[0]}\ldots \aamalg{\Delta[0]} \Delta[1] }_{n}\rightarrow \Delta[n]\]
    for $n\ge1$, and 
    
    \item the \emph{completeness cofibration}, given by either inclusion of the form
    \[\Delta[0]\to N\mathbb I.\]
\end{enumerate}
Complete Segal spaces, the fibrant objects in $\spsh{\Delta}_{(\infty,1)}$, are then precisely the injectively fibrant simplicial spaces that are local with respect to the maps of type (1) and (2).

\begin{rmk}
As briefly addressed in \cite[\textsection10]{rezkTheta}, in presence of the maps of type (1), for the purpose of the localization one could replace the map of type (2) as completeness acyclic cofibration with
\begin{enumerate}[leftmargin=*]
    \item[(2')] either inclusion of the form
    \[\Delta[0]\to\Delta[0]\aamalg{\Delta[1]}\Delta[3]\aamalg{\Delta[1]}\Delta[0],\]
where the right-hand side is the colimit of the diagram
\[\Delta[0]\xleftarrow{}\Delta[1]\xrightarrow{02}\Delta[3]\xleftarrow{13}\Delta[1]\xrightarrow{}\Delta[0].\]
\end{enumerate}

\end{rmk}

The homotopy theories of Segal categories and complete Segal spaces are equivalent.

\begin{thm}[{\cite[Thm 6.3]{bergner3models}}]
\label{secatcssqe}
The inclusion functor
from the category of Segal precategories to the category of simplicial spaces
induces a left Quillen equivalence
\[I\colon P\cat_{(\infty,1)}\to\spsh{\Delta}_{(\infty,1)}.\]
\end{thm}

\subsection{Complete Segal $\bT$-spaces as a model of $(\infty,2)$-categories}

We now recall the notion of complete Segal $\bT$-spaces, which give a model for $(\infty,2)$-categories.  

Let $\bT$ be Joyal's cell category.  For a precise account on how $\bT$ is defined we refer the reader to the original source \cite{JoyalDisks}, or to \cite[Def.\ 3.3]{BergerIterated} or \cite[\textsection 1.1]{rezkTheta} for an inductive approach; we give a brief review here.

Recall that $\bT$ is a full subcategory of $2\cat$ and that a generic object of $\bT$ is a $2$-category $[m|k_1, \ldots, k_m]$ generated by gluing horizontally the suspensions of $[k_i]$ for $i=1, \ldots, m$.  An example is the $2$-category $[3|2,0,1]$, which is generated by the following data:
\[
\begin{tikzcd}[row sep=3.2cm, column sep=2.2cm]
  x \arrow[r, bend left=50, "f", ""{name=U,inner sep=2pt,below}]
  \arrow[r, "g"{near end, xshift=0.2cm}, ""{name=D,inner sep=2pt},""{name=M,inner sep=2pt, below}]
  \arrow[r, bend right=50, "h"{below}, ""{name=DD,inner sep=2pt, xshift=0.05cm}]
  & y\arrow[r, "l"]
    & [-0.9cm]
    %& 
    z\arrow[r, bend left=50, "m", ""{name=U1,inner sep=2pt,below}]
  \arrow[r, "k"{near end, xshift=0.2cm}, ""{name=D1,inner sep=2pt, xshift=0.05cm},""{name=M1,inner sep=2pt, below}]
  & w.
    \arrow[Rightarrow, from=U, to=D, "\alpha"]
    \arrow[Rightarrow, from=M, to=DD, "\beta"{near start}]  
  \arrow[Rightarrow, from=U1, to=D1, "\gamma"]
\end{tikzcd}\]

\begin{defn}
A \emph{$\bT$-set} is a presheaf $A\colon \bT^{\op}\to\set$, and we denote the category of $\bT$-sets and natural transformations by $\psh{\bT}$.  Similarly, a \emph{$\bT$-space} is a simplicial presheaf $A\colon \bT^{\op}\to\sset$, and we denote the category of $\bT$-spaces by $\spsh{\bT}$.
\end{defn}

\begin{rmk}
The reader familiar with Rezk's paper \cite{rezkTheta} might observe that we are using the term ``$\bT$-space'' in a more general sense than he does.  His $\bT$-spaces satisfy additional Segal and completeness conditions that we discuss below; we further specify such objects by calling them ``complete Segal $\bT$-spaces''.
\end{rmk}

\begin{rmk}
The canonical inclusion $\set\hookrightarrow\sset$ of sets as discrete simplicial sets induces a canonical inclusion $\psh{\bT}\hookrightarrow\spsh{\bT}$, which is both a left and right adjoint. In particular, we often regard $\bT$-sets as discrete $\bT$-spaces without further specification.
\end{rmk}

\begin{notn}
For any $\theta$ in $\bT$, we denote by $\bT[\theta]$ the $\bT$-set represented by $\theta$.
\end{notn}

\begin{rmk}\label{BoxProduct}
As a special case of \cite[\textsection 3.1]{Ara}, given any $\bT$-set $A$ and any space $B$ one can consider the $\bT$-space $A\boxtimes B$, which is defined levelwise as the simplicial set
\[(A\boxtimes B)_{\theta}:=A_{\theta}\times B.\]
The construction extends to a bifunctor
\[\boxtimes\colon\psh{\bT}\times\sset\to\spsh{\bT}\]
that preserves colimits in each variable.
\end{rmk}

In preparation for a localization on the category $\spsh{\bT}$, we introduce the following class of maps. The reader may notice the analogy with the maps treated in \cref{segalspacesection}.
 
\begin{defn} \label{anodynemapstheta}
An \emph{elementary acyclic cofibration} is a map of discrete $\bT$-spaces of the following kinds.
\begin{enumerate}[leftmargin=*]
\item A \emph{vertical Segal acyclic cofibration} is given by, for some $k \ge0$, the canonical map
\[ \bT[1|1]\aamalg{\bT[1|0]}\dots\aamalg{\bT[1|0]}\bT[1|1]\hookrightarrow\bT[1|k].\]

\item A \emph{horizontal Segal acyclic cofibration} is given by, for some $m\ge0$ and $k_i\ge0$, where $0 \le i \le m$, the canonical map
\[ \bT[1|k_1]\aamalg{\bT[0]}\dots\aamalg{\bT[0]}\bT[1|k_m]\hookrightarrow\bT[m|k_1,\dots,k_m].\]

\item The \emph{horizontal completeness acyclic cofibration}
is either of the inclusions of the form
\[\bT[0]\to\bT[0]\aamalg{\bT[1]}\bT[3|0,0,0]\aamalg{\bT[1]}\bT[0],\]
where the right-hand side is the colimit of the diagram
\[\bT[0]\xleftarrow{}\bT[1]\xrightarrow{02}\bT[3|0,0,0]\xleftarrow{13}\bT[1] \xrightarrow{}\bT[0].\]
\item The \emph{vertical completeness acyclic cofibration} is the canonical map
\[\bT[1|0]\to\bT[1|0]\aamalg{\bT[1|1]}\bT[1|3]\aamalg{\bT[1|1]}\bT[1|0],\]
 induced by suspending the previous one.
\end{enumerate}
\end{defn}

We now describe two model structures on the category $\spsh{\bT}$, both established by Rezk \cite[\textsection 2.13, Prop.~11.5]{rezkTheta}.  Our description, in terms of the elementary acyclic cofibrations defined above, differs slightly from his, but is designed to facilitate some of our proofs in the next section.  We explain in \cref{RezkCompleteness} why the two approaches give the same model structures.

\begin{thm}
\label{RezkMS} The category $\spsh{\bT}$ of $\bT$-spaces supports the following two cofibrantly generated model structures:
\begin{itemize}[leftmargin=*]
\item the model structure $\spsh{\bT}_{i,(\infty,2)}$ obtained by taking the left Bousfield localization of the injective model structure $\spsh{\bT}_{inj}$ with respect to the set of elementary acyclic cofibrations from \cref{anodynemapstheta}; and
\item the model structure $\spsh{\bT}_{p,(\infty,2)}$ obtained by taking the left Bousfield localization of the projective model structure $\spsh{\bT}_{proj}$ with respect to the set of elementary acyclic cofibrations from \cref{anodynemapstheta}.
\end{itemize}
\end{thm}

Although the model structure $\spsh{\bT}_{i,(\infty,2)}$ is more common in the literature, for technical reasons that we discuss in \cref{DiscreteLCof}, in this paper we mostly work with the $\spsh{\bT}_{p,(\infty,2)}$. In this model structure,
\begin{itemize}[leftmargin=*]
\item the fibrant objects, which we call \emph{complete Segal $\bT$-spaces}, are precisely the projectively fibrant $\bT$-spaces that are local with respect to the elementary acyclic cofibrations from \cref{anodynemapstheta}; and
\item the cofibrations are precisely the projective cofibrations.
\end{itemize}

\begin{rmk} \label{ProjectiveGenerating}
Combining \cite[Thm 11.6.1, Def.\ 11.5.33, Def.\ 11.5.25]{Hirschhorn}, we obtain an explicit description of the generating cofibrations and generating acyclic cofibrations of $\spsh{\bT}_{p,(\infty,2)}$. 
More precisely:
\begin{enumerate}[leftmargin=*]
    \item a set of generating cofibrations for the projective model structure on $\spsh{\bT}$ is given by all maps of the form
    \[ \bT[\theta]\boxtimes\partial\Delta[\ell]\to\bT[\theta]\boxtimes\Delta[\ell]\quad\text{ for }\theta \in \Ob(\bT) \mbox{ and }\ell\geq 0;\]
    \item a set of generating acyclic cofibrations for the projective model structure on $\spsh{\bT}$ is given by all maps of the form
    \[\bT[\theta]\boxtimes\Lambda^k[\ell]\to\bT[\theta]\boxtimes\Delta[\ell]\quad\text{ for }\theta \in \Ob(\bT) \mbox{ and } 0\leq k \leq \ell.\]
  \end{enumerate}
\end{rmk}

The following equivalence between the two model structures can alternatively also be seen as a direct application of \cite[Theorem 3.3.20]{Hirschhorn}. 

\begin{thm}[{\cite[\textsection\textsection 2.5-2.13]{rezkTheta}}] \label{ProjInjEqui}
The identity functor defines a Quillen equivalence
\[\spsh{\bT}_{p,(\infty,2)}\rightleftarrows\spsh{\bT}_{i,(\infty,2)}.\]
\end{thm}

We want to consider the suspension of a simplicial space to a $\bT$-space.  Rezk uses the notation $V[1](X)$ in \cite[\textsection 4.4]{rezkTheta} for what we denote here by $\Sigma X$ to emphasize the analogy with similar constructions we have discussed.

\begin{defn}
The \emph{suspension} $\Sigma X$ of a simplicial space $X$ is the $\bT$-space obtained by applying the cocontinous functor $\Sigma\colon\spsh{\Delta}\to\spsh{\bT}_{*,*}$ defined on representable simplicial spaces as
\[\Sigma(\Delta[k]\boxtimes\Delta[\ell]):=\bT[1|k]\boxtimes\Delta[\ell].\]
This construction extends to a functor $\Sigma\colon\spsh{\Delta}\to\spsh{\bT}_{*,*}$ valued in bipointed  $\bT$-spaces.
\end{defn}

\begin{rmk}\label{RezkCompleteness}
In Rezk's original construction from \cite[\textsection 2.13, Prop.~11.5]{rezkTheta}, two model structures on $\spsh{\bT}$ are obtained by localizing the injective and projective model structure with respect to the set of maps of the following kinds:
\begin{enumerate}[label=(\arabic*'), leftmargin=*]
\item a family of maps that can be recognized to be precisely the family of vertical Segal acyclic cofibrations, using \cite[Prop.\ 11.7]{rezkTheta}; 

\item a family of maps that can be recognized to be precisely the family of horizontal Segal acyclic cofibration, using \cite[Prop.\ 11.7]{rezkTheta};
    
\item the unique map
     \[\Ntheta\mathbb I\to\bT[0]; \]
and

\item the map
 \[\Sigma\Ntheta\mathbb I\to\bT[1|0]\]
 obtained by suspending the map from (3').
\end{enumerate}

However, in presence of the maps of type (1) and (2), it is shown in \cite[\textsection 10]{rezkTheta} and also in \cite[\textsection 13]{BarwickSchommerPries} that for the purpose of the localization the maps of type (3) and (4) are equivalent to the maps of type (3') and (4'), respectively. It follows that, although presented differently, these two model structures in fact agree with the model structures $\spsh{\bT}_{i,(\infty,2)}$ and $\spsh{\bT}_{p,(\infty,2)}$ from \cref{RezkMS}.
\end{rmk}

\subsection{Complicial sets as a model of $(\infty,2)$-categories}

The first model of $(\infty,2)$-categories that we consider is based on the following structure, originally referred to as a \emph{simplicial set with hollowness} in \cite{StreetOrientedSimplexes} and later as a \emph{stratified simplicial set} in \cite{VerityComplicialAMS}.

\begin{defn}
A \emph{simplicial set with marking} is a simplicial set endowed with a subset of simplices of strictly positive dimensions that contain all degenerate simplices, called \emph{thin} or \emph{marked}. We denote by $m\sset$ the category of simplicial sets with marking and marking preserving simplicial maps.
\end{defn} 

We will consider a model structure on the category of simplicial sets with marking, in which the fibrant objects, called \emph{$2$-complicial sets}, provide a model for $(\infty,2)$-categories.
The idea is that, in a $2$-complicial set, the marked $k$-simplices are precisely the $k$-equivalences. We refer the reader to \cite{EmilyNotes} for further elaboration on this viewpoint.

\begin{rmk}
As discussed in \cite[Obs.~ 97]{VerityComplicialAMS}, the underlying simplicial set functor $\msset\to\sset$ fits into an adjoint triple
\[
\begin{tikzcd}[column sep=2cm]
\msset \arrow[r, ""{name=x1, above}, ""{name=x2, below, inner sep=1pt}]& \sset. \arrow[l, bend left=25, "(-)^{\sharp}"{ pos=0.47}, ""{name=x3, above,pos=0.47, inner sep=0pt}] \arrow[l, bend right=25, "(-)^{\flat}"{name=x4, above}, ""{name=x5, below, inner sep=1pt, pos=0.5}]
 \arrow[from=x1, to=x3, symbol=\dashv]
 \arrow[from=x5, to=x2, symbol=\dashv]
\end{tikzcd}
\]
For any simplicial set $X$, the left adjoint $X^\flat$ (sometimes also denoted simply by $X$) is obtained by marking only the degenerate simplices of $X$, and the right adjoint $X^{\sharp}$ is obtained by marking all simplices in positive dimensions.
\end{rmk}

\begin{rmk}
As described in detail in \cite[Obs.~109]{VerityComplicialAMS}, the category $m\sset$ of simplicial sets with marking is complete and cocomplete, with limits and colimits constructed as follows.
\begin{itemize}[leftmargin=*]
\item The underlying simplicial set of a limit $\lim_{i\in I}X_i$ of simplicial sets with marking is the limit of the corresponding underlying simplicial sets of $X_i$, and a simplex is marked in a limit of simplicial sets with marking $\lim_{i\in I}X_i$ if and only  if it is marked in each component $X_i$ for $i\in I$.

\item The underlying simplicial set of a colimit $\colim_{i\in I}X_i$ of simplicial sets with marking is the colimit of the corresponding underlying simplicial sets of $X_i$, and a simplex is marked in a colimit of simplicial sets with marking $\colim_{i\in I}X_i$ if and only if it admits a marked representative in $X_i$ for some $i\in I$.
\end{itemize}
\end{rmk}

The following model structure is one instance of the family of model structures constructed by Verity in \cite[Thm~100]{VerityComplicialI}, and is described in more detail in \cite[\textsection 3.3]{EmilyNotes}.

\begin{thm}[{\cite[Thm~1.25]{or}}] \label{modelstructureonstrat}
The category $\msset$ of simplicial sets with marking supports a cofibrantly generated cartesian closed model structure in which
\begin{itemize}[leftmargin=*]
\item the fibrant objects are the $2$-complicial sets, as recalled in \cite[Def.~1.21]{or}, and

\item the cofibrations are precisely the monomorphisms on underlying simplicial sets.
\end{itemize}
We denote this model structure by $\msset_{(\infty,2)}$.
\end{thm}

We warn the reader that the fibrant objects in this model structure have been given different names in the literature, and could perhaps more accurately be called
``$2$-trivial saturated weak complicial sets''. We have chosen to call them ``2-complicial sets'' for the sake of brevity; in what follows we do not make explicit use of their definition. We recall the key results we need, in particular about the weak equivalences in this model structure, in the remainder of this section. 

\begin{rmk} \label{SaturationWE}
Because of the way the model structure $\msset_{(\infty,2)}$ is constructed, if $\eqDelta$ denotes the $3$-simplex $\Delta[3]$ in which the non-degenerate marked $1$-simplices are precisely the one between the vertices $0$ and $2$ and the one between the vertices $1$ and $3$, and all simplices in dimension $2$ or higher are marked, the canonical map $\eqDelta\to\Delta[3]^{\sharp}$ is a weak equivalence. Indeed, the model structure $\msset_{(\infty,2)}$ is a Cisinski--Olschok model structure (in the sense of \cite{Olschok}) for which the map $\eqDelta\to\Delta[3]^{\sharp}$ is an anodyne extension.
\end{rmk}

\begin{lem}
\label{SharpLeftQuillen}
The functor
\[(-)^{\sharp}\colon\sset_{(\infty,0)}\to\msset_{(\infty,2)} \]
is a left Quillen functor, where $\sset_{(\infty,0)}$ denotes the Kan--Quillen model structure on the category $\sset$.
\end{lem}

\begin{proof}
The fact that the functor admits a right adjoint, often called the core functor, is discussed in \cite[Def.\ D.1.2]{RiehlVerityBook}. It is straightforward from its description that the functor $(-)^{\sharp}$ preserves cofibrations, and it is shown in \cite[Lem.~2.16]{or} that it also sends acyclic cofibrations of $\sset_{(\infty,0)}$ to weak equivalences of $\msset_{(\infty,2)}$. It follows that $(-)^{\sharp}$ defines indeed a left Quillen functor between the desired model categories.
\end{proof}

For $n=2$, the Street nerve was studied in detail by Duskin in \cite{duskin}, and can be described explicitly as follows.

\begin{defn}
The \emph{nerve} $N\cD$ of a $2$-category $\cD$ is the $3$-coskeletal simplicial set in which
\begin{enumerate}[leftmargin=*]
\item[(0)] a $0$-simplex consists of an object of $\cD$:
$$x;$$
    \item a $1$-simplex consists of a $1$-morphism of $\cD$:
     \[
%\left\{
\begin{tikzcd}
    x \arrow[rr, "a"{below}]&& y;
\end{tikzcd}
\]
    \item a $2$-simplex consists of a $2$-cell of $\cD$ of the form $c\Rightarrow b\circ a$:
    $$\begin{tikzcd}[baseline=(current  bounding  box.center)]
 & y \arrow[rd, "b"]  & \\%{\alpha_s}
    x \arrow[ru, "{a}"]
  \arrow[rr, "c"{below}, ""{name=D,inner sep=1pt}]
  && z;
  \arrow[Rightarrow, from=D, 
 to=1-2, shorten >= 0.1cm, shorten <= 0.1cm, "\varphi"]
\end{tikzcd}$$
\item a $3$-simplex consists of four $2$-cells of $\cD$ that satisfy the following relation:
\[
\simpfverArrow{d}{c.}{e}{a}{b}{f}{}{}{}\simpfverArrowcontinued{}{x}{y}{z}{w}
\]
\end{enumerate}
and in which the simplicial structure is as indicated in the pictures.
\end{defn}

\begin{defn}[{\cite[Ch.~10]{VerityComplicialAMS}}]
The \emph{Roberts--Street nerve} of a $2$-category $\cD$ is the simplicial set with marking $\NRS\cD$, defined by the following properties.
\begin{enumerate}[leftmargin=*]
    \item[(0)] The underlying simplicial set is the Duskin nerve $N\cD$. 
    
    \item[(1)] Only degenerate $1$-simplices are marked.
    
    \item[(2)] A $2$-simplex of $N\cD$ is marked in $\NRS\cD$ if and only if corresponding $2$-morphism $\varphi\colon c\Rightarrow b\circ a$ is an identity.
    \item[(3)] Any $m$-simplex of $N\cD$ for $m\ge3$ is marked in $\NRS\cD$.
\end{enumerate}
This\ construction extends to a functor $\NRS\colon2\cat\to\msset$.
 \end{defn}

The Roberts--Street nerve is a right adjoint functor, but, as proved by the second- and third-named authors, does not 
preserve fibrant objects on the model structures we want to consider.  However, it is
a \emph{homotopical} functor between model categories, in the sense that preserves weak equivalences. 

\begin{prop}[{\cite[Prop.\ 1.18]{ORfundamentalpushouts}}]  \label{NRSproperties}
The Roberts--Street nerve defines
a homotopical functor of model categories
\[\NRS\colon 2\cat\to \msset_{(\infty,2)}.\]
\end{prop}

The following two technical results essentially tell us that horizontal and vertical composition of $2$-cells can be encoded via Segal-type maps that are acyclic cofibrations in the model structure for 2-complicial sets.

\begin{thm}[{\cite[Cor.~2.10]{ORfundamentalpushouts}}]\label{StratHorizontalSegality} 
For any $m\ge0$ and $k_i\ge0$ for $i=1,\dots,m$ there is a canonical map of simplicial sets with marking
\[\NRS [1|k_1]\aamalg{\NRS [0]}\dots\aamalg{\NRS [0]}\NRS [1|k_m]\hookrightarrow \NRS [m|k_1,\dots,k_m] \]
that is an acyclic cofibration, and in particular a weak equivalence, in $\msset_{(\infty,2)}$.
\end{thm}

\begin{thm}[{\cite[Cor.~2.11]{ORfundamentalpushouts}}]\label{StratVerticalSegality}
For any $k\ge0$ there is a canonical map of simplicial sets with marking
\[\NRS [1|1]\aamalg{\NRS [1|0]}\dots\aamalg{\NRS [1|0]}\NRS [1|1]\hookrightarrow \NRS [1|k]\]
that is an acyclic cofibration, and in particular a weak equivalence, in $\msset_{(\infty,2)}$.
\end{thm}

An important construction in this paper is the suspension of a simplicial set with marking.  We conclude this section with the definition and some key results about it.
 
\begin{defn}[{\cite[Def.~2.6]{ORfundamentalpushouts}}]
The \emph{suspension} $\Sigma X$ of a simplicial set with marking $X$ is the simplicial set with marking defined as follows.
\begin{itemize}[leftmargin=*]
    \item It has precisely two $0$-simplices,
    which we call $x_{\bot}$ and $x_{\top}$, respectively.
    \item The set of $m$-simplices for $m>0$ is given by all $k$-simplices of $X$ for $0\le k\le m-1$ as well as the $m$-fold degeneracies of the two $0$-simplices $x_{\bot}$ and $x_{\top}$, namely
    \[(\Sigma X)_m\cong\{s_0^mx_{\bot}\}\amalg X_{m-1} \amalg \ldots \amalg X_0 \amalg \{s_0^mx_{\top}\}.\]
    \item The simplicial structure can be read off from \cite[Def.~2.6]{ORfundamentalpushouts}.
    \item The set of non-degenerate $m$-simplices for $m>0$ is given by the non-degenerate $(m-1)$-simplices of $X$.
\item A non-degenerate $m$-simplex $\sigma$ is marked in $\Sigma X$ if and only if it is marked as an $(m-1)$-simplex of $X$.
\end{itemize} 
This construction extends to a functor $\Sigma\colon\msset\to\msset_{*,*}$ valued in bipointed marked simplicial sets.
\end{defn}

We now recall that this functor can be upgraded to a left Quillen functor of model categories. Recall from \cite{HirschhornOvercategories} that, given any cofibrantly model category $\cM$, there is a model structure on the category $\cM_{*,*}$ of bipointed objects in $\cM$, in which cofibrations, fibrations, and weak equivalences are created in $\cM$.

\begin{lem}[{\cite[Lem.~2.7]{ORfundamentalpushouts}}]
\label{suspensionhomotopical}
Regarding $\Sigma X$ as a simplicial set with marking bipointed on $x_{\bot}$ and $x_{\top}$, the marked suspension defines a left Quillen functor
\[\Sigma\colon\msset_{(\infty,2)}\to (\msset_{(\infty,2)})_{*,*}.\]
In particular, it is homotopical and it respects connected colimits as a functor $\Sigma\colon m\sset\to\msset$.
\end{lem}

Finally, we recall that the suspension of a marked simplicial set is homotopically compatible with the Roberts--Street nerve, as one would expect.  

\begin{thm}[{\cite[Thm~2.9]{ORfundamentalpushouts}}] \label{StratNerveSuspension}
For any $1$-category $\cD$ there is a canonical natural map
\[\Sigma \NRS \cD\to \NRS \Sigma\cD\]
that is a weak equivalence in $\msset_{(\infty,2)}$.
\end{thm}
\section{The comparison of models of $(\infty,2)$-categories}
\label{ComparisonStructure}

In this section, we set up our explicit comparison between the two models for $(\infty,2)$-categories that we are considering.  We first establish the desired Quillen pair of functors between the  unlocalized model structure on the category of $\bT$-spaces and the model structure on simplicial sets with marking, then show that is is still a Quillen pair after localizing the former model category.  We then show that it is a Quillen equivalence, deferring some steps in the proof to later sections.

\subsection{The Quillen pair before localizing}

Let us begin by defining the functor that we use to make our comparison.

\begin{const}
The functor $\bT\times\Delta\subset\spsh{\bT}\to\msset$ given by 
\[(\theta,[\ell])\mapsto(\bT\times\Delta)[\theta,\ell]=\bT[\theta]\boxtimes\Delta[\ell]\mapsto \NRS\theta\times\Delta[\ell]^\sharp.\]
induces an adjunction
\[L\colon\spsh{\bT}\rightleftarrows\msset\colon R.\]
\end{const}

Roughly speaking, for any $\bT$-space $W$, the simplicial set with marking $LW$ is obtained by gluing together a copy of the Roberts--Street nerve of $\theta$, for any $\theta$ in $\bT$ that maps to $W$.  While describing this gluing explicitly is complicated, it is essentially specified by the definition of left Kan extension. 

We now show that these adjoint functors define Quillen pair on unlocalized model categories.

\begin{prop} \label{Linjective}
The adjunction
\[L\colon\spsh{\bT}_{proj}\rightleftarrows\msset_{(\infty,2)}\colon R\]
is a Quillen pair.
\end{prop}

\begin{proof}
We want to show that the functor $L$ preserves cofibrations and acyclic cofibrations.  From \cref{ProjectiveGenerating} we know that
\begin{enumerate}[leftmargin=*]
    \item a set of  generating cofibrations for the projective model structure on $\spsh{\bT}$ is given by all maps of the form
    \[\bT[\theta]\boxtimes\partial\Delta[\ell]\to\bT[\theta]\boxtimes\Delta[\ell]\quad\text{ for }\theta \in \bT \mbox{ and }\ell\geq 0;\]
    and 
    
    \item a set of generating acyclic cofibrations for the projective model structure on $\spsh{\bT}$ is given by all maps of the form
    \[\bT[\theta]\boxtimes\Lambda^k[\ell]\to\bT[\theta]\boxtimes\Delta[\ell]\quad\text{ for }\theta \in \bT \mbox{ and } 0\leq k \leq \ell.\]
  \end{enumerate}
Using the facts that $(-)^{\sharp}$ commutes with colimits, which is a consequence of \cref{SharpLeftQuillen}, and that the box product $\boxtimes$ preserves colimits in each variable, which was recalled in \cref{BoxProduct}, we see that
\begin{enumerate}[leftmargin=*]
\item the image of the generating cofibration via $L$ is the map
 \[\NRS\theta\times\partial\Delta[\ell]^\sharp\to\NRS\theta\times\Delta[\ell]^\sharp\quad\text{ for }\theta \in \bT \mbox{ and }\ell\geq 0;\]
 
\item the image of the generating acyclic cofibration via $L$ is the map
\[\NRS\theta\times\Lambda^k[\ell]^\sharp\to\NRS\theta\times\Delta[\ell]^\sharp\quad\text{ for }\theta \in \bT \mbox{ and } 0\leq k \leq \ell.\]
    \end{enumerate}
    Since the model structure $\msset_{(\infty,2)}$ is cartesian closed by \cref{modelstructureonstrat} and $(-)^{\sharp}$ is a left Quillen functor by \cref{SharpLeftQuillen}, we conclude that
    \begin{enumerate}[leftmargin=*]
    \item the map $\NRS\theta\times\partial\Delta[\ell]^\sharp\to\NRS\theta\times\Delta[\ell]^\sharp$ is a cofibration and
    \item the map $\NRS\theta\times\Lambda^k[\ell]^\sharp\to\NRS\theta\times\Delta[\ell]^\sharp$ is an acyclic cofibration
    \end{enumerate}
    It follows that $L$ preserves cofibrations and acyclic cofibrations, so it is a left Quillen functor, as desired.
\qedhere
\end{proof}

\begin{rmk} \label{DiscreteLCof}
One might wonder, in contrast with much of the literature on the subject, why we have chosen to use the projective, rather than the injective, model structure on $\spsh{\bT}$.  However, it is not clear whether the functor 
\[L\colon\spsh{\bT}_{inj}\to\msset_{(\infty,2)}\]
is a left Quillen functor, since we do not know whether it preserves cofibrations. More precisely, it is unclear whether $L$ sends the injective cofibration
\[\partial\bT[3|1,0,1]\to\bT[3|1,0,1]\]
to a cofibration of $\msset_{(\infty,2)}$.
\end{rmk}

\subsection{The Quillen pair after localizing}

We now show that we still have a Quillen pair after localizing the projective model structure on $\spsh{\bT}$.

\begin{thm} \label{LQuillenPair}
The adjunction
\[L\colon\spsh{\bT}_{p,(\infty,2)}\rightleftarrows\msset_{(\infty,2)}\colon R\]
is a Quillen pair.
\end{thm}

Since cofibrations are unchanged by localization, it suffices to prove that $L$ preserves acyclic cofibrations.  We do so by proving that $L$ preserves all elementary acyclic cofibrations, in the following sequence of propositions.

\begin{prop}
The functor $L$ sends the vertical Segal acyclic cofibrations 
\[ \bT[1|1]\aamalg{\bT[1|0]}\dots\aamalg{\bT[1|0]}\bT[1|1]\hookrightarrow\bT[1|k]\text{ for }k\ge0 \]
 from \cref{anodynemapstheta} to weak equivalences in $\msset_{(\infty,2)}$.
\end{prop}

\begin{proof}
The functor $L$ sends the elementary acyclic cofibration
 \[\bT[1|1]\aamalg{\bT[1|0]}\dots\aamalg{\bT[1|0]}\bT[1|1]\hookrightarrow\bT[1|k]\]
to the canonical inclusion
\[\NRS [1|1]\aamalg{\NRS [1|0]}\dots\aamalg{\NRS [1|0]}\NRS [1|1]\hookrightarrow \NRS [1|k]\]
which is an acyclic cofibration by \cref{StratVerticalSegality}.
\end{proof}

\begin{prop}
The functor $L$ sends the horizontal Segal acyclic cofibrations 
\[ \bT[1|k_1]\aamalg{\bT[0]}\dots\aamalg{\bT[0]}\bT[1|k_m] \hookrightarrow \bT[m|k_1,\dots,k_m]\text{ for }m\ge0\text{ and }k_i\ge0 \]
from \cref{anodynemapstheta} to weak equivalences in $\msset_{(\infty,2)}$.
\end{prop}

\begin{proof}
The functor $L$ sends the elementary acyclic cofibration
\[\bT[1|k_1]\aamalg{\bT[0]}\dots\aamalg{\bT[0]}\bT[1|k_m]\hookrightarrow\bT[m|k_1,\dots,k_m]\]
to the canonical inclusion
\[ \NRS[1|k_1]\aamalg{\NRS[0]}\dots\aamalg{\NRS[0]}\NRS[1|k_m]\hookrightarrow \NRS[m|k_1,\dots,k_m]\]
which is an acyclic cofibration by \cref{StratHorizontalSegality}.
\end{proof}

\begin{prop} \label{1Complete}
The functor $L$ sends the horizontal completeness acyclic cofibration 
 \[\bT[0]\to\bT[0]\aamalg{\bT[1|0]}\bT[3|0,0,0]\aamalg{\bT[1|0]}\bT[0]\]
from \cref{anodynemapstheta} to a weak equivalence of $\msset_{(\infty,2)}$.
\end{prop}

To prove this proposition, we need the following preliminary lemma.

\begin{lem} \label{QuotientSaturation}
The unique map
\[\Delta[0]\aamalg{\Delta[1]}\NRS[3]\aamalg{\Delta[1]}\Delta[0]\to\Delta[0]\]
is a weak equivalence in $\msset_{(\infty,2)}$.
\end{lem}

\begin{proof}
We observe that this map fits into a commutative diagram of simplicial sets with marking
\[
\begin{tikzcd}
\Delta[0]\aamalg{\Delta[1]}\NRS[3]\aamalg{\Delta[1]}\Delta[0] \arrow[r] & \Delta[0]\\
\Delta[1]_t\aamalg{\Delta[1]}\NRS[3]\aamalg{\Delta[1]}\Delta[1]_t \cong\eqDelta\arrow[r]\arrow[u]& \Delta[3]^{\sharp}\arrow[u].
\end{tikzcd}
\]
In this diagram, we observe that:
\begin{itemize}[leftmargin=*]
    \item the bottom horizontal map is an acyclic cofibration by \cref{SaturationWE};
    
    \item the left vertical map is a map between (homotopy) pushouts induced by the identity of $\NRS[3]$ and two copies of the weak equivalence $\Delta[1]_t\to\Delta[0]$;
    and
    
    \item the right vertical map is a weak equivalence since $(-)^{\sharp}$ preserves weak equivalences by \cref{SharpLeftQuillen}.
\end{itemize}
It follows by two-out-of-three that the top horizontal map is a weak equivalence, as desired.
\end{proof}

We can now use this lemma to prove \cref{1Complete}.

\begin{proof}[Proof of \cref{1Complete}]
The functor $L$ sends the map
\[\bT[0]\to\bT[0]\aamalg{\bT[1|0]}\bT[3|0,0,0]\aamalg{\bT[1|0]}\bT[0]\]
to a map
\[\Delta[0]\to\Delta[0]\aamalg{\Delta[1|0]}\NRS[3|0,0,0]\aamalg{\Delta[1|0]}\Delta[0]\]
which we want to show is a weak equivalence. However, we can conclude this fact by the two-out-of-three property, since we know from \cref{QuotientSaturation} that the unique map
\[\Delta[0]\aamalg{\Delta[1]}\NRS[3]\aamalg{\Delta[1]}\Delta[0]\to\Delta[0]\]
is a weak equivalence in $\msset_{(\infty,2)}$.
\end{proof}

To complete the proof of \cref{LQuillenPair}, it remains to show that $L$ preserves one more acyclic cofibration.

\begin{prop}\label{2Complete}
The functor $L$ sends the vertical completeness acyclic cofibration from \cref{anodynemapstheta}
\[\bT[1|0]\to\bT[1|0]\aamalg{\bT[1|1]}\bT[1|3]\aamalg{\bT[1|1]}\bT[1|0]\]
to a weak equivalence of $\msset_{(\infty,2)}$.
\end{prop}

\begin{proof}
The functor $L$ sends the map
\[\bT[1|0]\to\bT[1|0]\aamalg{\bT[1|1]}\bT[1|3]\aamalg{\bT[1|1]}\bT[1|0]\]
to a map
\[\NRS [1|0]\to \NRS [1|0]\aamalg{\NRS [1|1]}\NRS[1|3]\aamalg{\NRS [1|1]}\NRS [1|0]\]
which we want to show is a weak equivalence. By the two-out-of-three property, it suffices to show that the map 
\[\NRS [1|0]\aamalg{\NRS [1|1]}\NRS[1|3]\aamalg{\NRS [1|1]}\NRS [1|0]\to\NRS [1|0], \]
induced by the unique map $[1|3]\to[1|0]$ in $\bT$ that is bijective on objects, is a weak equivalence in $\msset_{(\infty,2)}$. This map can be rewritten in terms of suspensions of $1$-categories, as in \cref{suspension}, as
\[\NRS \Sigma[0]\aamalg{\NRS \Sigma[1]}\NRS\Sigma[3]\aamalg{\NRS \Sigma[1]}\NRS \Sigma[0]\to\NRS\Sigma[0].\]
By \cref{StratNerveSuspension}, this map fits into a commutative diagram of simplicial sets with marking
\[
\begin{tikzcd}
\NRS \Sigma[0]\aamalg{\NRS \Sigma[1]}\NRS\Sigma[3]\aamalg{\NRS \Sigma[1]}\NRS \Sigma[0]\arrow[r]& \NRS\Sigma[0]\\
\Sigma\NRS [0]\aamalg{\Sigma\NRS [1]}\Sigma\NRS[3]\aamalg{\Sigma\NRS [1]}\Sigma\NRS [0]\arrow[r]\arrow[u]&\arrow[u,"\cong"]\Sigma\NRS[0],
\end{tikzcd}
\]
in which the two vertical maps are weak equivalences.  Note that for the left-hand map we are using the fact that these pushouts are actually homotopy pushouts. In particular, by the two-out-of-three property, to prove the theorem it is enough to prove that the bottom map is a weak equivalence.
Using the fact that suspension commutes with pushouts by \cref{suspensionhomotopical}, this map can be rewritten as
\[\Sigma\left(\Delta[0]\aamalg{\Delta[1]}\NRS[3]\aamalg{\Delta[1]}\Delta[0]\right)\to\Sigma\Delta[0],\]
namely the suspension of the map
\[\Delta[0]\aamalg{\Delta[1]}\NRS[3]\aamalg{\Delta[1]}\Delta[0]\to\Delta[0]\]
which was shown in \cref{QuotientSaturation} to be a weak equivalence. Since suspension is a left Quillen functor by \cref{suspensionhomotopical}, we are done.
\end{proof}

\subsection{The Quillen equivalence}

It remains to show that this Quillen pair is in fact a Quillen equivalence.  Our proof, however, is not done directly via the definition, but instead uses some machinery due to Barwick and Schommer-Pries in \cite{BarwickSchommerPries} that we now briefly recall.

The first thing we need to consider is their criterion for when a model category is a ``model of $(\infty,2)$-categories''.  We begin with the some notation.

\begin{notn} 
Given a model category $\cM$, we denote by $\cM_{\infty}$ the underlying $(\infty,1)$-category of $\cM$.  While we do not need one here, for explicit (different but equivalent) constructions of $\cM_{\infty}$ in the model of quasi-categories, we refer the reader to \cite{HinichDKLoc} or \cite[\textsection A.2]{htt}.
\end{notn}

\begin{notn}
Given a Quillen pair $F\colon\cM\rightleftarrows\cM'\colon G$ between model categories $\cM$ and $\cM'$, we denote by
\[F_{\infty}\colon\cM_{\infty}\rightleftarrows\cM'_{\infty}\colon G_{\infty}\]
the adjunction that $(F,G)$ induces at the level of underlying $(\infty,1)$-categories.

On objects, the value of $F_{\infty}$ on any object of $\cM$ can be computed up to equivalence in $\cM'_{\infty}$ by applying $F$ to any cofibrant replacement of the given object. Similarly the value of $G_{\infty}$ on any object of $\cM'$ can be computed up to equivalence in $\cM_{\infty}$ by applying $G$ to any fibrant replacement of the given object.
Moreover, if $(F,G)$ is a Quillen equivalence, the induced adjunction $(F_{\infty},G_{\infty})$ is an equivalence of $(\infty,1)$-categories.  For more details on how to obtain this adjunction of $(\infty,1)$-categories in the model of quasi-categories we refer the reader to \cite[Prop.~1.5.1]{HinichDKLoc}. 
\end{notn}

\begin{defn}[Barwick--Schommer-Pries]
A model category $\cM$ is a \emph{model of $(\infty,2)$-categories} if the underlying $(\infty,1)$-category is equivalent to the colossal model $\cK$ from \cite[\textsection8]{BarwickSchommerPries}, namely if there exists an equivalence of $(\infty,1)$-categories
\[\cM_{\infty}\simeq \cK.\]
\end{defn}

The \emph{colossal model}, which we denote by $\cK$ here for simplicity, is constructed as an $(\infty,1)$-category in \cite[\textsection 8]{BarwickSchommerPries}.
As we discuss in the appendix, with standard techniques one can also present the colossal model as the underlying $(\infty,1)$-category of a model category. More precisely, we show as \cref{modelcategorycolossal} that it is the underlying $(\infty,1)$-category $(\spsh{\Upsilon_2}_{(\infty,2)})_{\infty}$ of a model category $\spsh{\Upsilon_2}_{(\infty,2)}$.
 
In any case, for the main purpose of this paper the arguments are packaged in a way that no explicit construction for the colossal model is needed.

\begin{thm}
The model categories $\spsh{\bT}_{p,(\infty,2)}$ and $\spsh{\bT}_{i,(\infty,2)}$ are models of $(\infty,2)$-categories.
\end{thm}

\begin{proof}
The fact that $\spsh{\bT}_{i,(\infty,2)}$ is a model of $(\infty,2)$-categories is discussed in \cite[\textsection 13]{BarwickSchommerPries} and there is an equivalence
\[\left(\spsh{\bT}_{p,(\infty,2)}\right)_{\infty}\simeq\left(\spsh{\bT}_{i,(\infty,2)}\right)_{\infty}\]
induced by the Quillen equivalence from \cref{ProjInjEqui}.
\end{proof}

\begin{thm}
The model category $\msset_{(\infty,2)}$ is a model of $(\infty,2)$-categories.
\end{thm}

\begin{proof}
An equivalence of $(\infty,1)$-categories
\[\left(\vphantom{\spsh{\bT}_{i,(\infty,2)}}\msset_{(\infty,2)}\right)_{\infty}\simeq\left(\spsh{\bT}_{i,(\infty,2)}\right)_{\infty}\]
can be obtained combining several equivalences of $(\infty,1)$-categories induced by Quillen equivalences due to \cite{JT, lurieGoodwillie, br1, br2, GHL}, as we recalled in \eqref{DiagramIntro}.
\end{proof}

Next, we recall the definition of a $j$-cell in a model of $(\infty,2)$-categories.

\begin{defn}\label{CellDef} 
Let $\cM$ be a model category that is a model for $(\infty,2)$-categories. An object of $\cM$ is a representative of the \emph{$j$-cell} for $j=0,1,2$ if it corresponds to the $j$-cell of the colossal model through any equivalence of $(\infty,1)$-categories $\cM_{\infty}\simeq(\spsh{\Upsilon_2}_{(\infty,2)})_{\infty}$.
\end{defn}

For completeness, the definitions of the \emph{$0$-, $1$- and $2$-cells} in the colossal model are recalled in \cref{AxiomaticAppendix}, but will not be needed explicitly.

\begin{rmk}\label{CellsWellDef}
The object in $\cM$ that represents the $j$-cell is unique up to equivalence in $\cM_{\infty}$ and also up to isomorphism in $\Ho\cM$, the homotopy category of $\cM$. The definition makes sense in particular because any auto-equivalence of $(\spsh{\Upsilon_2}_{(\infty,2)})_{\infty}$ preserves $j$-cells for $j=0,1,2$, as shown in \cite[Thm 7.3]{BarwickSchommerPries}.
\end{rmk}

The following statements describe $j$-cells in $\spsh{\bT}_{p,(\infty,2)}$ and $\msset_{(\infty,2)}$.

\begin{restatable}{propRestate}{RestateCellsRezkProj}
\label{CellsRezkproj}
In $\spsh{\bT}_{p,(\infty,2)}$ the object $\bT[C_j]$
is a representative of the $j$-cell for $j=0,1,2$.
\end{restatable}

\begin{restatable}{propRestate}{RestateStratCells}
\label{StratCells}
In $\msset_{(\infty,2)}$ the object $\NRS C_j$
is a representative of the $j$-cell for $j=0,1,2$.
\end{restatable}

Although the two statements are not surprising, the argument to identify cells in $\msset_{(\infty,2)}$ requires significant work and makes use of many external results. We therefore postpone both proofs to \cref{Sec:RecognizingCells}.

Finally, the following theorem is the key ingredient to prove that the functor $L$ is a Quillen equivalence.

\begin{thm}[{\cite[Prop.~15.10]{BarwickSchommerPries}}] \label{BSPCriterion}
Let $\cM$ and $\cN$ be model categories that are models for $(\infty,2)$-categories, and $L\colon\cM\rightleftarrows\cN\colon R$ a Quillen pair between them. Then the Quillen pair $(L,R)$ is a Quillen equivalence if and only if the derived functor of $L$ sends $j$-cells to $j$-cells for $j=0,1,2$.
\end{thm}

Once the proofs of \cref{CellsRezkproj,StratCells} are provided in \cref{Sec:RecognizingCells}, we can then apply \cref{BSPCriterion} to the Quillen pair from \cref{LQuillenPair} to conclude the desired Quillen equivalence.

\begin{thm}
The adjunction
\[L\colon\spsh{\bT}_{p,(\infty,2)}\rightleftarrows\msset_{(\infty,2)}\colon R\]
is a Quillen equivalence, and in particular induces an equivalence of $(\infty,1)$-categories
\[L_{\infty}\colon\left(\spsh{\bT}_{p,(\infty,2)}\right)_{\infty}\rightleftarrows\left(\vphantom{\spsh{\bT}_{i,(\infty,2)}}\msset_{(\infty,2)}\right)_{\infty}\colon R_{\infty}.\]
\end{thm}

\section{Recognizing cells in models of $(\infty,2)$-categories}
\label{Sec:RecognizingCells}

The goal of this section is to identify the $j$-cells in $\spsh{\bT}_{p,(\infty,2)}$, and most importantly the $j$-cells in $\msset_{(\infty,2)}$, as defined in \cref{CellDef}.
The structure of the argument involves the identification of the $j$-cells in several established model categories that are models of $(\infty,2)$-categories.

In the following picture, we display the equivalences used to identify the cells in the marked simplicial sets, and the propositions displayed show how the cells behave under the corresponding equivalence:
\begin{equation}
    \begin{tikzpicture}[baseline=-2.7cm]
    \def\r{4.5cm}
\node[outer sep=0.2cm] (A1) at (0,-0.4*\r) {$\spsh{\bT}_{i,(\infty,2)}$};
\node[outer sep=0.2cm] (A2) at (0,0) {$\spsh{\bT}_{p,(\infty,2)}$};
\node (A3) at (0,-0.8*\r) {$\spsh{(\Delta\times\Delta)}_{(\infty,2)}$};
\node (A4) at (0,-1.2*\r) {$P\cat(\spsh{\Delta})_{(\infty,2)}$};
\node (A5) at (0.9*\r,-1.2*\r) {$\vcat{\spsh{\Delta}_{(\infty,1)}}$};
\node (B1) at (1.8*\r, -1.2*\r) {$\vcat{\psh{\Delta}_{(\infty,1)}}$}; 
\node (C1) at (1.8*\r, -0.8*\r) {$\vcat{\sset^+_{(\infty,1)}}$}; 
\node (C2) at (1.8*\r, -0.4*\r) {$\sset^{sc}_{(\infty,2)}$};
\node (D1) at (1.8*\r,0) {$\msset_{(\infty,2)}$.};

 \path (A1) -- node(comp1)[sloped,transform shape,inner sep=0pt,xscale=2.5] {$\simeq$} 
        node[font=\tiny, right of=comp1]{Prop.\ \ref{CellsRezkproj}}(A2);
 \path (A1) -- node(comp2)[sloped,transform shape,inner sep=0pt,xscale=2.5, rotate=180] {$\simeq$} 
        node[font=\tiny, right of=comp2]{Prop.\ \ref{CellTwofoldSegal}}(A3);
\path (A3) -- node(comp3)[sloped,transform shape,inner sep=0pt,xscale=2.5, rotate=180] {$\simeq$} 
        node[font=\tiny, right of=comp3]{Prop.\ \ref{CellsTwoSegalCat}}(A4);   
%%%%
\path (A4) -- node(comp4)[sloped,transform shape,inner sep=0pt,xscale=2.5, outer sep=0pt] {$\simeq$} 
        node[font=\tiny, below of=comp4, yshift=0.4cm]{Prop.\ \ref{CellsCSSenriched}}(A5);   
\path (A5) -- node(comp5)[sloped,transform shape,inner sep=0pt,xscale=2.5, outer sep=0pt] {$\simeq$} 
        node[font=\tiny, below of=comp4, yshift=0.4cm]{Prop.\ \ref{CellsJoyalenriched}}(B1); 
\path (B1) -- node(comp6)[sloped,transform shape,inner sep=0pt,xscale=2.5, rotate=180] {$\simeq$} 
        node[font=\tiny, left of=comp6]{Prop.\ \ref{CellsMarkedenriched}}(C1);
\path (C1) -- node(comp7)[sloped,transform shape,inner sep=0pt,xscale=2.5, rotate=180] {$\simeq$}
        node[font=\tiny, left of=comp7]{Prop.\ \ref{CellsScaled}}(C2);
\path (C2) -- node(comp8)[sloped,transform shape,inner sep=0pt,xscale=2.5, rotate=180] {$\simeq$} 
        node[font=\tiny, left of=comp7]{Prop.\ \ref{StratCells}}(D1);
    \end{tikzpicture}
\end{equation}

While it is not possible to make this section completely self-contained, we have included precise references for all relevant constructions and definitions.

\begin{lem} \label{CellRecognition}
Suppose that a functor $F\colon\cM\to\cM'$ is a left (respectively, right) Quillen equivalence between models of $(\infty,2)$-categories, and an object $X$ is cofibrant (respectively, fibrant) in $\cM$. Then $X$ is a $j$-cell in $\cM$ for some $0\leq j \leq 2$ if and only if $F(X)$ is a $j$-cell in $\cM'$.
\end{lem}

\begin{proof}
Consider the induced functor $F_{\infty} \colon \cM_{\infty} \to \cM'_{\infty}$, which is an equivalence of $(\infty,1)$-categories.  It follows that, for any $j=0,1,2$ and $j$-cell $X_j$ of $\cM$, the object $F_{\infty}(X_j)$ is a cell in $\cM'_{\infty}$, either by direct verification, or by appealing to \cref{BSPCriterion}.  Now, an object $X$ is a $j$-cell in $\cM$ if and only if there is an isomorphism $X\cong X_j$ in $\Ho\cM$. Again using the fact that $F_{\infty}$ is an equivalence, this statement is equivalent to saying that there is an isomorphism $F_{\infty}(X) \cong F_{\infty}(X_j)$ in $\Ho(\cM')$ . But the existence of such an isomorphism is equivalent to having $F_{\infty}(X)$ a $j$-cell of $\cM'_{\infty}$ because $F_{\infty}(X_j)$ is one.  Since $F(X)$ computes $F_{\infty}(X)$, the result follows.
\end{proof}

\subsection{Recognizing cells in $\bT$-models of $(\infty,2)$-categories}

We now begin the work of identifying $j$-cells in different models for $(\infty,2)$-categories.  We begin with the $j$-cells in $\spsh{\bT}_{i,(\infty,2)}$, which have been identified by Barwick and Schommer-Pries.

\begin{prop}[{\cite[\textsection 13]{BarwickSchommerPries}}]\label{CellRezkInj}
In $\spsh{\bT}_{i,(\infty,2)}$ the object $\bT[C_j]$
is a representative of the $j$-cell for $j=0,1,2$.
\end{prop}

We can prove \cref{CellsRezkproj}, that identifies the cells in $\spsh{\bT}_{p,(\infty,2)}$.

\begin{proof}[Proof of \cref{CellsRezkproj}] 
We consider the identity functor on $\spsh{\bT}$, which is by \cref{ProjInjEqui} a left Quillen equivalence \[\id\colon\spsh{\bT}_{p,(\infty,2)}\to\spsh{\bT}_{i,(\infty,2)}.\]
For $j=0,1,2$, by \cref{CellRezkInj} we know that $\bT[C_j]$ is a $j$-cell in $\spsh{\bT}_{i,(\infty,2)}$.  Moreover, the object $\bT[C_j]$ is projectively cofibrant by \cite[Prop.~11.6.2]{Hirschhorn}, since it is representable. It then follows from \cref{CellRecognition} that $\bT[C_j]$ is a $j$-cell in $\spsh{\bT}_{p,(\infty,2)}$.
\end{proof}

%%%
%%%
%%%

%%%
%%%
%%%

\subsection{Recognizing cells in multisimplicial models of $(\infty,2)$-categories}

We now turn to identifying $j$-cells in multisimplicial models of $(\infty, 2)$-categories.  Because we have not yet considered these models in this paper, we describe them briefly as we go.

\begin{thm}[{\cite[Ch.~2]{BarwickThesis}}]
The category $\spsh{(\Delta\times\Delta)}$ of bisimplicial spaces admits a model structure in which
\begin{itemize}[leftmargin=*]
    \item the fibrant objects are the complete Segal objects in complete Segal spaces; and
    
    \item the cofibrations are the monomorphisms, and in particular every object is cofibrant.
\end{itemize}
We denote this model structure by $\spsh{(\Delta\times\Delta)}_{(\infty,2)}$.
\end{thm}

The idea behind complete Segal objects in complete Segal spaces is that we apply similar Segal and completeness conditions to functors $\Deltaop \to \spsh{\Delta}$, where the target category is equipped with the complete Segal space model structure.  Thus, the Segal and completeness maps are now equivalences in this model structure, rather than equivalences of simplicial sets. For more details on the definition of \emph{complete Segal objects in complete Segal spaces}, see \cite[Ch.~2]{BarwickThesis}, \cite[Def.~5.3]{br2}, or \cite{lurieGoodwillie}.

There is an explicit equivalence of model categories between this model structure and the one for complete Segal $\bT$-spaces.  See also \cite{BarwickSchommerPries} for a different proof of this equivalence.

\begin{thm}[{\cite[Cor.~7.1]{br2}}]
The functor $d\colon\Delta\times\Delta\to\bT$ given by $[m,k]\mapsto[m|k,\dots,k]$ induces a left Quillen equivalence
\[d^*\colon\spsh{\bT}_{i,(\infty,2)}\to\spsh{(\Delta\times\Delta)}_{(\infty,2)}.\]
\end{thm}

In particular, $\spsh{(\Delta\times\Delta)}_{(\infty,2)}$ is a model for $(\infty,2)$-categories.  We now characterize the $j$-cells in this model.

\begin{prop}\label{CellTwofoldSegal}
In $\spsh{(\Delta\times\Delta)}_{(\infty,2)}$ the object $d^*\bT[C_j]$
is a representative of the $j$-cell for $j=0,1,2$.
\end{prop}

\begin{proof}
We consider the functor $d^*$, which is a left Quillen equivalence 
\[ d^*\colon\spsh{\bT}_{i,(\infty,2)}\to\spsh{(\Delta\times\Delta)}_{(\infty,2)}.\]
By \cref{CellRezkInj}, we know that for $j=0,1,2$, the object $\bT[C_j]$ is a $j$-cell in $\spsh{\bT}_{i,(\infty,2)}$.  Moreover, every object is cofibrant in $\spsh{\bT}_{i,(\infty,2)}$. It follows from \cref{CellRecognition} that $d^*\bT[C_j]$ is a $j$-cell in $\spsh{(\Delta\times\Delta)}_{(\infty,2)}$.
\end{proof}

%%%
%%%
%%%

We can now generalize the notion of Segal precategory to this context; in analogy with the notion of complete Segal objects described above, we can define \emph{Segal precategory objects} in complete Segal spaces, given by functors $X \colon \Deltaop \rightarrow \spsh{\Delta}$ such that $X_0$ is a discrete object and the Segal maps are weak equivalences in the complete Segal space model structure.  See \cite[\textsection 6]{br1} for more details.

Let us briefly describe the comparison with complete Segal objects, which is analogous to \cref{secatcssqe}.  We denote by $P\cat(\psh{\Delta})$ the full subcategory of $\spsh{(\Delta\times\Delta)}$ given by \emph{Segal precategory objects} in simplicial spaces, namely those bisimplicial spaces $X \colon \Deltaop \rightarrow \spsh{\Delta}$ for which $X_0$ is discrete, and we denote by $I\colon P\cat(\psh{\Delta})\to\spsh{(\Delta\times\Delta)}$ the inclusion functor. 

\begin{thm}[{\cite[\textsection 6]{br1}}]
The category $P\cat(\spsh{\Delta})$ of precategories in simplicial spaces admits a model structure in which
\begin{itemize}[leftmargin=*]
    \item the fibrant objects are the Segal category objects;  and 
    
    \item the cofibrations are the monomorphisms, and in particular every object is cofibrant.
\end{itemize}
We denote this model structure by $P\cat(\spsh{\Delta})_{(\infty,2)}$.
\end{thm}

In particular, $P\cat(\spsh{\Delta})_{(\infty,2)}$ is a model for $(\infty,2)$-categories.

\begin{thm}[{\cite[Prop.~9.5, Thm~9.6]{br2}}]
The natural inclusion functor from \cite[\textsection9]{br2}
induces a left Quillen equivalence
\[I\colon P\cat(\spsh{\Delta})_{(\infty,2)}\to\spsh{(\Delta\times\Delta)}_{(\infty,2)}.\]
\end{thm}

In particular, $P\cat(\spsh{\Delta})_{(\infty,2)}$ is a model for $(\infty,2)$-categories.

For each $j=0,1,2$, the bisimplicial space $d^*\bT[C_j]$, a priori an object of $\spsh{(\Delta\times\Delta)}$, is actually a precategory, so it can be regarded as an object of $P\cat(\spsh{\Delta})$.

\begin{prop}\label{CellsTwoSegalCat}
In $P\cat(\spsh{\Delta})_{(\infty,2)}$ the object $d^*\bT[C_j]$
is a representative of the $j$-cell for $j=0,1,2$.
\end{prop}

\begin{proof}
We consider the inclusion functor, which is a left Quillen equivalence
\[I\colon P\cat(\spsh{\Delta})_{(\infty,2)}\to\spsh{(\Delta\times\Delta)}_{(\infty,2)}.\]
For $j=0,1,2$, by \cref{CellTwofoldSegal} we know that $I(d^*\bT[C_i])$ is a $j$-cell in $\spsh{(\Delta \times \Delta)}_{(\infty,2)}$. Moreover, every object is cofibrant in $P\cat(\spsh{\Delta})_{(\infty,2)}$. It follows from \cref{CellRecognition} that $d^*\bT[C_i]$ is a $j$-cell in $P\cat(\spsh{\Delta})_{(\infty,2)}$.
\end{proof}

%%%
%%%
%%%

%%%
%%%
%%%

\subsection{Recognizing cells in enriched models of $(\infty,2)$-categories}

We now turn to recognizing cells in models that are given by enriched categories.  Many model structures on enriched categories can be obtained by the following general result of Lurie.

\begin{thm}[{\cite[Thm A.3.2.24]{htt}}] \label{EnrichedCatMS}
Let $\cV$ be an excellent monoidal model category, in the sense of \cite[Def.\ A.3.2.16]{htt}. The category of small
categories enriched over $\cV$ admits a model structure in which
\begin{itemize}[leftmargin=*]
    \item the fibrant objects are the \emph{locally fibrant categories}, i.e., the enriched categories $C$ so that for any pair of objects $c,c'$ in $C$, the mapping object $C(c,c')$ is fibrant in $\cV$; and 
    
    \item the weak equivalences, which are described in \cite[Def.\ A.3.2.1]{htt} and \cite{LawsonEnriched}, are
    enriched functors $F\colon C \to D$ so that:
    \begin{enumerate}[leftmargin=*]
        \item for every pair of objects $c,c'$ of $C$, the map induced by $F$ of mapping objects 
        \[
        F_{c,d}\colon\Map_{C}(c,c') \to \Map_{D}(Fc,Fc'),
        \]
        is a weak equivalence in $\cV$, and
        \item the functor induced by $F$ on (underlying categories) of
        $\Ho\cV$-categories is essentially surjective.
    \end{enumerate}
    \item the cofibrations are those described in \cite[Prop.~A.3.2.4]{htt}.
\end{itemize}
We denote this model structure by $\vcat{\cV}$.
\end{thm}

To give an idea, the technical condition for a combinatorial monoidal model category to be \emph{excellent} requests a closure property for cofibrations and weak equivalences, additionally to compatibility of the model structure with the monoidal structure. Lurie's original definition also requires a further condition, known as ``invertibility hypothesis'', which was shown to follow from the other conditions by Lawson as \cite[Thm~0.1]{LawsonEnriched}.

We specialize this construction to the following situations. 
\begin{itemize}[leftmargin=*]
\item Let $\cV=\cat$ be the canonical model structure on the category $\cat$ of small categories from \cite{RezkCat}, which is seen to be excellent using the fact that the nerve functor creates weak equivalences and commutes with filtered colimits.  We then obtain precisely the model category $\vcat{\cat}=2\cat$ as discussed in \cite[Ex.\ 1.8]{BergerMoerdijkEnriched}.

\item Let $\cV=\sset_{(\infty,1)}$ be the Joyal model structure on the category $\sset$ of simplicial sets from \cite[Thm~6.12]{JoyalVolumeII}, which is excellent by \cite[Ex.\ A.3.2.23]{htt}.  We then obtain the model category $\vcat{\sset_{(\infty,1)}}$.

\item Let $\cV=\spsh{\Delta}_{(\infty,1)}$ being the Rezk model structure from \cref{Rezkmodelstructure}
%\cite[Thm~7.2]{rezk}
on the category $\spsh{\Delta}$ of simplicial spaces, which is discussed to be excellent \cite[Thm 3.11]{br1}.  We then obtain the model category $\vcat{\spsh{\Delta}_{(\infty,1)}}$.

\item Let $\cV=\sset^+_{(\infty,1)}$ be the Lurie model structure on the category $\sset^+$ of marked simplicial sets from \cite[Prop.~3.1.3.7]{htt}, which is excellent by \cite[Ex.\ A.3.2.22]{htt}.  We then obtain the model category $\vcat{\sset^+_{(\infty,1)}}$.
\end{itemize}

We now turn to an explicit Quillen equivalence between one of these enriched models and one of the models we have already discussed.

\begin{thm}[{\cite[7.1-7.6]{br1}}]
The enriched nerve functor from \cite[Def.~7.3]{br1}, obtained by regarding a bisimplicial category as a simplicial object in simplicial spaces, defines a right Quillen equivalence
\[R\colon\vcat{\spsh{\Delta}_{(\infty,1)}}\to P\cat(\spsh{\Delta})_{(\infty,2)}.\]
\end{thm}

In particular, $\vcat{\spsh{\Delta}_{(\infty,1)}}$ is a model for $(\infty,2)$-categories.
 
Now, we would like to identify the $j$-cells in the model structure $\vcat{\spsh{\Delta}_{(\infty,1)}}$, for which we make use of the discrete nerve functor $N^{\mathrm{disc}} \colon \cat \to \spsh{\Delta}$ considered in \cite{rezk}.  Since preserves products (being a right adjoint functor), it induces a functor $N^{\mathrm{disc}}_* \colon\vcat{\cat}\to\vcat{\spsh{\Delta}}$, given by applying $N^{\mathrm{disc}}$ to each mapping category.

\begin{prop}\label{CellsCSSenriched}
In $\vcat{\spsh{\Delta}_{(\infty,1)}}$ the object $N^{\mathrm{disc}}_*C_j$
is a representative of the $j$-cell for $j=0,1,2$.
\end{prop}

Before proving this proposition, we establish two lemmas that tell us more about the structure of these discrete nerves.

\begin{lem}\label{FibrancyCSSenrichedCells}
For any $j=0,1,2$, the $\spsh{\Delta}$-enriched category $N^{\mathrm{disc}}_{*}\bT[C_j]$ is fibrant in $\vcat{\spsh{\Delta}_{(\infty,1)}}$.
\end{lem}

\begin{proof}
For $j=0,1,2$, all hom-categories of $C_j$ are of the form $\varnothing=[-1]$, $[0]$, or $[1]$. So all hom-bisimplicial sets of $N^{\mathrm{disc}}_{*}\bT[C_j]$ are of the form $N^{\mathrm{disc}}[-1]$, $N^{\mathrm{disc}}[0]$ or $N^{\mathrm{disc}}[1]$, which are all complete Segal spaces, namely fibrant in $\spsh{\Delta}_{(\infty,1)}$, since the categories $[-1]$, $[0]$, and $[1]$ do not have any non-trivial isomorphisms.
\end{proof}

\begin{lem} \label{TwoPrecatNerves}
For any $\theta$ in $\bT$, there is an isomorphism of precategories
\[R(N^{\mathrm{disc}}_*\theta)\cong d^*\bT[\theta].\]
\end{lem}

\begin{proof}
For $i,j,k\ge0$, we first compute the set $\left(RN^{\mathrm{disc}}_*\theta\right)_{[i],[j],[k]}$.
If $\cD$ is a bisimplicial category with object set $\cD_0$, and $\cD_1$ denotes the bisimplicial space \[\cD_1=\coprod_{a,b \in \cD_0}\Map_{\cD}(a,b),\]
by definition of $R$ (as given in \cite[Def.~7.3]{br1}) for any $i\ge0$ there is an isomorphism of bisimplicial sets
\[\left(R\cD\right)_{[i]}\cong \underbrace{\cD_1 \underset{\cD_0}{\times}\cD_1 \underset{\cD_0}{\times} \ldots \underset{\cD_0}{\times} \cD_1}_{i}
,\]
that is natural in $i$. When specializing to the case $\cD=N^{\mathrm{disc}}_*\theta$, we obtain a natural isomorphism
\[\left(RN^{\mathrm{disc}}_*\theta\right)_{[i]}\cong \underbrace{(N^{\mathrm{disc}}_*\theta)_1 \underset{(N^{\mathrm{disc}}_*\theta)_0}{\times}(N^{\mathrm{disc}}_*\theta)_1 \underset{(N^{\mathrm{disc}}_*\theta)_0}{\times} \ldots \underset{(N^{\mathrm{disc}}_*\theta)_0}{\times} (N^{\mathrm{disc}}_*\theta)_1}_{i}.\]
In particular, if $\theta_0$ denotes the set of objects of $\theta$ and $\theta_1$ denotes the category
\[\theta_1:=\coprod_{a,b \in \theta_0}\Map_{\theta}(a,b),\]
for any $j,k\ge0$ we have a bijection
\[\left(RN^{\mathrm{disc}}_*\theta\right)_{[i],[j],[k]}\cong \underbrace{N_j\theta_1\underset{\theta_0}{\times} N_j\theta_1 \underset{\theta_0}{\times}\ldots \underset{\theta_0}{\times} N_j\theta_1}_{i}
,\]
that is natural in $i,j,k$.
 
Next, for $i,j,k\ge0$, we compute the set $(d^*\bT[\theta])_{[i],[j],[k]}$.
By definition of $d^*$, and using the fact that $\bT$ is a full subcategory of $\twocat$, we have bijections
\[\begin{array}{rcl}
\left(d^*\bT[\theta]\right)_{[i],[j],[k]}&\cong& \Hom_{\Theta_2}([i|\underbrace{j,j,\ldots,j}_{i}], \theta) \\
&\cong& \Hom_{\twocat}([i|\underbrace{j,j,\ldots,j}_{i}], \theta) \\
&\cong &\Hom_{\twocat}(\underbrace{[1|j]\aamalg{[0]}[1|j] \aamalg{[0]} \ldots \aamalg{[0]} [1|j]}_{i}, \theta) \\
&\cong &\underbrace{\Hom_{\twocat}([1|j], \theta)\underset{\Hom_{\twocat}([0], \theta)}{\times} \ldots \underset{\Hom_{\twocat}([0], \theta)}{\times} \Hom_{\twocat}([1|j], \theta)}_i,\\
&\cong &\underbrace{\Hom_{\twocat}([1|j], \theta)\underset{\theta_0}{\times} \ldots \underset{\theta_0}{\times} \Hom_{\twocat}([1|j], \theta)}_i,\\
\end{array}\]
that are natural in $i,j,k$.

Finally, we show that there is a bijection
\[ \Hom_{\twocat}([1|j], \theta)\cong N_j\theta_1\]
that is natural in $j$, from which the lemma follows.
To do so, we observe that there are natural bijections 
    \[\begin{array}{rcl}
    \Hom_{\twocat}([1|j], \theta)&\cong &\coprod_{a,b \in \theta_0} \Hom_{\twocat_{*,*}}([1|j], (\theta, a,b))\\
    &\cong& \coprod_{a,b \in \theta_0}\Hom_{\cat}([j], \Map_{\theta}(a,b))\\
     &\cong& \Hom_{\cat}([j], \coprod_{a,b \in \theta_0}\Map_{\theta}(a,b))\\
          &\cong& \Hom_{\cat}([j], \theta_1)\cong N_j\theta_1,
    \end{array}\]
 as desired.
\end{proof}

\begin{proof}[Proof of \cref{CellsCSSenriched}]
Consider the right Quillen equivalence
\[R\colon\vcat{\spsh{\Delta}_{(\infty,1)}}\to P\cat(\spsh{\Delta})_{(\infty,2)}.\]
By \cref{CellTwofoldSegal} and \cref{TwoPrecatNerves}, we know that for any $j=0,1,2$, the object $d^*\bT[C_j]\cong R(N^{\mathrm{disc}}_*C_j)$ is a $j$-cell in $P\cat(\spsh{\Delta})_{(\infty,2)}$.
Moreover, by \cref{FibrancyCSSenrichedCells} the object $N^{\mathrm{disc}}_*C_j$ is fibrant in $\vcat{\spsh{\Delta}_{(\infty,1)}}$. It follows from \cref{CellRecognition} that $N^{\mathrm{disc}}_*C_j$ is a $j$-cell in $\vcat{\spsh{\Delta}_{(\infty,1)}}$, as desired.
\end{proof}

We now compare the model structure for categories enriched in complete Segal spaces to the model structure for categories enriched in quasi-categories.

\begin{thm}
\label{changeofenrichment}
The functor induced by taking $(-)_0$ on each hom-simplicial space defines a right Quillen equivalence
\[\vcat{\spsh{\Delta}_{(\infty,1)}}\to \vcat{\psh{\Delta}_{(\infty,1)}}.\]
\end{thm}

In particular, $\vcat{\psh{\Delta}_{(\infty,1)}}$ is a model for $(\infty,2)$-categories.

\begin{proof}
The functor $p\colon\Delta\times \Delta\to\Delta$, defined by $[m,n]\mapsto[m]$, induces an adjoint triple
\[
\begin{tikzcd}[column sep=2cm]
\psh{\Delta} \arrow[r, ""{name=x1, above}, ""{name=x2, below, inner sep=1pt}, "p^*"{near end}]& \psh{(\Delta\times\Delta)}=\spsh{\Delta}. \arrow[l, bend left=25, "p_*"{ pos=0.47}, ""{name=x3, above,pos=0.68, inner sep=0pt}] \arrow[l, bend right=25, "p_!"{name=x4, above}, ""{name=x5, below, inner sep=1pt, pos=0.685}]
 \arrow[from=x1, to=x3, symbol=\dashv]
 \arrow[from=x5, to=x2, symbol=\dashv]
\end{tikzcd}
\]
where $p^*$ is given by precomposition with $p$, while $p_!$ and $p_*$ are the left and right Kan extensions along $p$, respectively. In particular, the functor $p^*$ is (strong) monoidal with respect to cartesian product because it is a right adjoint. Moreover, it is shown as \cite[Thm~4.11]{JT} that the adjunction
\[p^*\colon\psh{\Delta}_{(\infty,1)}\rightleftarrows\spsh{\Delta}_{(\infty,1)}\colon p_*\]
is a Quillen equivalence.
One can then apply \cite[Rmk~A.3.2.6]{htt} to obtain the desired Quillen equivalence, observing that $p_*$ is the functor $(-)_0$.
\end{proof}

We can now use this equivalence to identify the $j$-cells in $\vcat{\psh{\Delta}_{(\infty,1)}}$.

\begin{prop} \label{CellsJoyalenriched}
In $\vcat{\psh{\Delta}_{(\infty,1)}}$ the object $N_*C_j$
is a representative of the $j$-cell for $j=0,1,2$.
\end{prop}

\begin{proof}
Consider the right Quillen equivalence from  \cref{changeofenrichment}
\[\vcat{\spsh{\Delta}_{(\infty,1)}}\to\vcat{\psh{\Delta}_{(\infty,1)}}.\] 
By \cref{CellsCSSenriched,FibrancyCSSenrichedCells} we know that for each $j=0,1,2$, the object $N^{\mathrm{disc}}_*C_j$ is a $j$-cell in $\vcat{\spsh{\Delta}_{(\infty,1)}}$ and is fibrant. It follows from \cref{CellRecognition} that $N_*C_j$, the image of $N_*^{\mathrm{disc}}C_j$ under the above Quillen equivalence, is a $j$-cell in $\vcat{\psh{\Delta}_{(\infty,1)}}$.
\end{proof}

We now make a similar comparison between categories enriched in quasi-categories and categories enriched in marked simplicial sets.

\begin{thm}
The functor induced by taking the underlying simplicial set $U$ on each mapping object defines a right Quillen equivalence
\[U_*\colon\vcat{\sset^+_{(\infty,1)}}\to \vcat{\sset_{(\infty,1)}}.\]
\end{thm}

In particular, $\vcat{\sset^+_{(\infty,1)}}$ is a model for $(\infty,2)$-categories.

\begin{proof}
The desired right Quillen equivalence is an instance of \cite[Rmk~A.3.2.6]{htt} applied to the right Quillen equivalence
\[U\colon \sset^+_{(\infty,1)}\to\psh{\Delta}_{(\infty,1)}\]
from \cite[Thm~3.1.5.1]{htt}.
\end{proof}

Once again, our goal is to identify the $j$-cells in this model structure.  To do so, consider the flat nerve functor $N^{\flat}\colon\cat\to\sset^+$, obtained by regarding the nerve of a category in which the marked $1$-simplices are precisely those corresponding to identity morphisms in the category. One can check that the functor $N^{\flat}$ preserves finite cartesian products, from which we obtain an induced functor $N^{\flat}_*\colon\vcat{\cat}\to\vcat{\sset^+}$, given by applying $N^{\flat}$ to each mapping category.

\begin{prop} \label{CellsMarkedenriched}
In $\vcat{\sset^+_{(\infty,1)}}$the object $N^\flat_*C_j$
is a representative of the $j$-cell for $j=0,1,2$.
\end{prop}

We begin with a lemma establishing that these objects are fibrant.

\begin{lem} \label{flatnervecellsfibrant}
For $j=0,1,2$, the object $N_*^{\flat}C_j$ is fibrant in $\vcat{\sset^+_{(\infty,1)}}$.
\end{lem}

\begin{proof}
For $j=0,1,2$, all hom-marked simplicial sets of $N_*^{\flat}C_j$ are of the form $N^\flat[-1]$, $N^\flat[0]$ or $N^\flat[1]$, which are naturally marked quasi-categories, and therefore fibrant in $\vcat{\sset^+_{(\infty,1)}}$, since the categories $[-1]$, $[0]$ and $[1]$ have no non-trivial isomorphisms.
\end{proof}

\begin{proof}[Proof of \cref{CellsMarkedenriched}]
We consider the right Quillen equivalence  \[U_*\colon\vcat{\sset^+_{(\infty,1)}}\to\vcat{\sset_{(\infty,1)}}.\]
By \cref{CellsJoyalenriched} we know that for each $j=0,1,2$, the object $U_*N_*^\flat C_j\cong N_*C_j$ is a $j$-cell in $\vcat{\sset_{(\infty,1)}}$, and $N^\flat_*C_j$ is fibrant in $\vcat{\sset^+_{(\infty,1)}}$ by \cref{flatnervecellsfibrant}. It follows from \cref{CellRecognition} that $N^\flat_*C_j$ is a $j$-cell in $\vcat{\sset^+_{(\infty,1)}}$.
\end{proof}
%%%
%%%
%%%

\subsection{Recognizing cells in simplicial models of $(\infty,2)$-categories}

Finally, we want to identify the $j$-cells in the model of marked simplicial sets.  To aid in doing so, we look first at the related model of scaled simplicial sets. A \emph{scaled simplicial set} is a simplicial set with a collection of marked $2$-simplices including degenerate $2$-simplices.

\begin{thm}[{\cite[Thm~4.2.7]{lurieGoodwillie}}]
The category $\sset^{sc}$ of scaled simplicial sets admits a model structure in which
\begin{itemize}[leftmargin=*]
    \item the fibrant objects are the $\infty$-bicategories from \cite[Def.~4.2.8]{lurieGoodwillie}, and
    \item the cofibrations are the monomorphisms (and in particular every object is cofibrant).
\end{itemize}
We denote this model structure by $\sset^{sc}_{(\infty,2)}$.
\end{thm}

Lurie enhances the classical homotopy coherent nerve functor $\mathfrak N\colon \vcat{\sset}\to\sset$ to the context of scaled simplicial sets by taking into account the marking, obtaining a scaled homotopy coherent nerve functor $\mathfrak N\colon \vcat{\sset^{+}}\to\sset^{sc}$.

\begin{thm}[{\cite[Thm 0.0.3]{lurieGoodwillie}}] 
The scaled homotopy coherent nerve functor from \cite[Def.\ 3.1.10]{lurieGoodwillie} defines a right Quillen equivalence
\[\mathfrak{N}^{sc}\colon\vcat{{\sset^+}_{(\infty,1)}}\to\sset^{sc}_{(\infty,2)}.\]
\end{thm}

In particular, $\sset^{sc}_{(\infty,2)}$ is a model for $(\infty,2)$-categories.  We now describe the $j$-cells in this model structure.

The description of the $j$-cells in this model structure makes use of a similar scaled nerve construction, in the form of a functor $N^{sc}\colon2\cat\to\sset^{sc}$, as described in \cite[Def.~8.1]{GHL}. Given any $2$-category $\cD$, the scaled nerve $N^{sc}\cD$ is given by the Duskin nerve of $\cD$ together with the marking of all $2$-simplices arising from $2$-isomorphisms.

\begin{prop} \label{CellsScaled}
In $\sset^{sc}_{(\infty,2)}$ the object $N^{sc}C_j$
is a representative of the $j$-cell for $j=0,1,2$.
\end{prop}

\begin{proof}
As a preliminary observation, we mention that there is an isomorphism of scaled simplicial sets
\[N^{sc}C_j\cong \mathfrak{N}^{sc}N_*^\flat C_j\]
for $j=0, 1,2$. This fact can be deduced combining \cite[Def.\ 8.1]{GHL} together with \cite[Prop.~8.2]{GHL}.

Consider now the right Quillen equivalence \[\mathfrak{N}^{sc}\colon\vcat{{\sset^+}_{(\infty,1)}}\to{\sset}^{sc}_{(\infty,2)}.\]
By \cref{CellsMarkedenriched} we know that for each $j=0,1,2$, the object $N_*^\flat C_j$ is a $j$-cell in $\vcat{\sset^+_{(\infty,1)}}$.  We have also proved that $N^\flat_*C_j$ is fibrant in $\vcat{\sset^+_{(\infty,1)}}$ in \cref{flatnervecellsfibrant}. It follows from \cref{CellRecognition} that $N^{sc}C_j\cong \mathfrak{N}^{sc}N^\flat_*C_j$ is a $j$-cell in $\sset^{sc}_{(\infty,2)}$.
\end{proof}

%%%
%%%
%%%

Finally, we can compare the models of scaled simplicial sets and marked simplicial sets.

\begin{thm}[{\cite[Thm 7.7]{GHL}}]
The forgetful functor defines a right Quillen equivalence
\[U\colon\msset_{(\infty,2)}\to{\sset}^{sc}_{(\infty,2)}.\]
\end{thm}

We can now prove \cref{StratCells}, which characterizes the cells in $\msset_{(\infty,2)}$.

\begin{proof}[Proof of \cref{StratCells}]
We consider the right Quillen equivalence \[U\colon\msset_{(\infty,2)}\to{\sset}^{sc}_{(\infty,2)}.\]
By \cite[Def.\ 8.1]{GHL} we know that $U\NRS C_j\cong N^{sc}C_j$ for each $j=0,1,2$. Moreover, we know that $N^{sc}C_j$ is a $j$-cell in ${\sset}^{sc}_{(\infty,2)}$ by \cref{CellsMarkedenriched} and that it is fibrant in $\msset_{(\infty,2)}$ by \cite[Thm~5.1(1)]{or}.  It follows from \cref{CellRecognition} that $\NRS C_j$ is a $j$-cell in $\msset_{(\infty,2)}$.
\end{proof}

\appendix \label{AxiomaticAppendix}

\stepcounter{section}

\section*{Appendix: The colossal model of $(\infty,2)$-categories}
\renewcommand{\thethm}{\Alph{section}.\arabic{thm}}
\renewcommand{\thenotn}{\Alph{section}.\arabic{thm}}
\renewcommand{\thedefn}{\Alph{section}.\arabic{thm}}
\renewcommand{\theprop}{\Alph{section}.\arabic{thm}}
\renewcommand{\thermk}{\Alph{section}.\arabic{thm}}

In this section, we give a model categorical variant of the colossal model by Barwick--Schommer-Pries.

In order to recall the original definition of the colossal model, we fix the following notations. We denote by $\Upsilon_2$ Barwick--Schommer-Pries' indexing category for the colossal model, namely the full subcategory of $2\cat$ from \cite[Def.~6.2]{BarwickSchommerPries}.
In particular, ${\left(\Upsilon_2^{\op}\right)_{\infty}}$ is the $(\infty,1)$-category obtained by regarding the category $\Upsilon_2^{\op}$ as an $(\infty,1)$-category.
We denote by $\cS_{\infty}$ the $(\infty,1)$-category of spaces, namely $\cS_{\infty}=(\sset_{(\infty,0)})_{\infty}$.

\begin{defn}[\cite{BarwickSchommerPries}]
The \emph{colossal model} is the $(\infty,1)$-category \[\cL_{\infty}\left({\cS}_{\infty}^{{\left(\Upsilon_2^{\op}\right)_{\infty}}}\right),\]
obtained by localizing the presheaf $(\infty,1)$-category ${\cS}_{\infty}^{(\Upsilon_2^{\op})_{\infty}}$ at the set of maps from \cite[Notn~8.3]{BarwickSchommerPries}.
\end{defn}

From the definition, we see that the colossal model is obtained by considering the $(\infty,1)$-category of spaces, taking a presheaf $(\infty,1)$-category valued in it, and the localizing. By contrast, one could instead present the same $(\infty,1)$-category by considering the Quillen model structure, which presents the $(\infty,1)$-category of spaces, taking the injective model structure on a presheaf category of functors valued in the Quillen model structure, and then a left Bousfield localization of it.
More precisely, one can consider the following model structure.

\begin{prop} 
The category $\spsh{\Upsilon_2}$ of $\Upsilon_2$-spaces supports a cofibrantly generated model structure obtained by taking the left Bousfield localization of the injective model structure $\spsh{\Upsilon_2}_{inj}$ with respect to the set of elementary acyclic cofibrations from \cite[Notn~6.5]{BarwickSchommerPries}.
We denote this model structure by $\spsh{\Upsilon_2}_{(\infty,2)}$.
\end{prop}

We want to prove that this model structure does present the colossal model, in the sense of the following theorem.

\begin{thm}
\label{modelcategorycolossal}
There is an equivalence of $(\infty,1)$-categories
\[\cL_{\infty}\left(({\sset_{(\infty,0)}})_{\infty}^{(\Upsilon_2^{\op})_{\infty}}\right)\simeq\left(\spsh{\Upsilon_2}_{(\infty,2)}\right)_{\infty}\]
\end{thm}
 
The proof is an application of the following result, which guarantees that one can build localizations of presheaf categories either as model categories or directly as $(\infty,1)$-categories.
 
\begin{prop} \label{LocFunInfty}
Let $\cA$ be a category, $\cM$ a left proper combinatorial simplicial model category, and $\Lambda$ a set of maps in $\cM^{\cA}$. There is an equivalence of $(\infty,1)$-categories 
\[\cL_{\infty}(\cM_{\infty}^{\cA_\infty})\simeq(\cL(\cM^{\cA}_{inj}))_{\infty},\]
where $\cL(\cM^{\cA}_{inj})$ denotes the Bousfield localization of the injective model structure $\cM^{\cA}_{inj}$ at $\Lambda$, and $\cL_{\infty}(\cM_{\infty}^{\cA_\infty})$ denotes the localization of the $\infty$-category $\cM_{\infty}^{\cA_\infty}$ at $\Lambda$.
\end{prop}

The proof of the proposition requires the following two ingredients.

\begin{thm}[{\cite[Proof of Prop.~A.3.7.8]{htt}}]\label{LocalizeInfty}
Let $\cN$ a left proper combinatorial simplicial model category, and $\Lambda$ a set of maps in $\cN$.
There is an equivalence of $(\infty,1)$-categories
\[\cL_{\infty}\cN_{\infty}\simeq(\cL\cN)_{\infty},\]
where $\cL\cN$ denotes the Bousfield localization of the model structure $\cN$ at $\Lambda$, and $\cL_{\infty}\cN_{\infty}$ denotes the localization of the $(\infty,1)$-category $\cN_{\infty}$ at $\Lambda_{\infty}$.
\end{thm}

\begin{thm}[{\cite[Prop.~4.2.4.4]{htt}}]\label{FunInfty}
Let $\cA$ be a category and $\cM$ a combinatorial simplicial model category.
There is an equivalence of $(\infty,1)$-categories
\[\cM_{\infty}^{\cA_\infty}\simeq(\cM^{\cA}_{inj})_{\infty},\]
where $\cM^{\cA}_{inj}$ denotes the injective model structure on $\cM^{\cA}$.
\end{thm}

We can now prove the proposition.

\begin{proof}[Proof of \cref{LocFunInfty}]
Combining \cref{LocalizeInfty,FunInfty}, we obtain an equivalence of $(\infty,1)$-categories
\[\cL_{\infty}(\cM_{\infty}^{\cA_{\infty}})\simeq \cL_{\infty}(\cM^{\cA})_{\infty}\simeq (\cL(\cM^{\cA}))_{\infty},
\]
as desired.
\end{proof}

We can now prove the theorem.

\begin{proof}[Proof of \cref{modelcategorycolossal}]
Applying \cref{LocFunInfty} with $\cM=\cS=\sset_{(\infty,0)}$, $\cA=\Upsilon_2^{\op}$ and $\Lambda=S$, the set of maps in ${\sset}^{\Upsilon_2^{\op}}$ from \cite[Notn~8.3]{BarwickSchommerPries}, we obtain the equivalence of $(\infty,1)$-categories
\[\cL_{\infty}\left(\cS_{\infty}^{(\Upsilon_2^{\op})_{\infty}}\right)=\cL_{\infty}\left(({\sset_{(\infty,0)})}_{\infty}^{(\Upsilon_2^{\op})_{\infty}}\right)\simeq\left(\cL\left(\sset_{(\infty,0)}^{\Upsilon_2^{\op}}\right)\right)_{\infty}=\left(\spsh{\Upsilon_2}_{(\infty,2)}\right)_{\infty}\]
as desired.
\end{proof} 

\bibliographystyle{amsalpha}
\bibliography{ref}

\end{document}

%% file: PreambleSep2020.tex
\usepackage{amsmath}
\usepackage{latexsym,amsfonts,amssymb,mathrsfs}
\usepackage{mathdots}
\usepackage{color}
\usepackage[all,2cell,color]{xy}
\usepackage{enumitem} %allows e.g. to change labels in enumerations
\usepackage{bbm}%allows to make mathbb-like numbers
\usepackage{tikz}
\usepackage{tikz-cd}
\usepackage[T1]{fontenc}
 \usepackage[utf8]{inputenc}
 \usepackage{blkarray} %package for matrices with labels
 \usepackage{mathtools} %package for alignment in matrices
 \usepackage{longtable}%allows to make reasonable tables, in particular with pagebreaks
 \usepackage{multicol}%allows to write in two columns
 
 \usepackage{thmtools}%these two packages allow to restate a theorem
\usepackage{thm-restate}
 
 \usepackage{hyperref}
\usepackage[capitalise]{cleveref}

\usetikzlibrary{positioning} %should allow to draw more precise arrows
\usetikzlibrary{decorations.pathreplacing}%need this to draw nice curly brackets
\usetikzlibrary{arrows, decorations.markings} %need this for nice arrows in tikz pictures, markings in particular for double arrows
\usetikzlibrary{calc} %allows to do calculations in tikz

\UseAllTwocells
\usepackage[
textwidth=3cm,
textsize=small,
colorinlistoftodos]
{todonotes}

\tikzset{%
    symbol/.style={%
        draw=none,
        every to/.append style={%
            edge node={node [sloped, allow upside down, auto=false]{$#1$}}}
    }
}
\tikzset{dot/.style={circle,draw,fill,inner sep=1pt},
    arrowinline/.style={-stealth,thick,shorten <=2pt,shorten >=2pt}}

\theoremstyle{plain}   % This is the default, anyway
\newtheorem{thm}{Theorem}[section] % numbered theorem
\makeatletter\let\c@thm\c@thm\makeatother

\makeatletter\let\c@cor\c@thm\makeatother
\newtheorem{lem}{Lemma}[section]
\makeatletter\let\c@lem\c@thm\makeatother
\newtheorem{prop}{Proposition}[section]
\makeatletter\let\c@prop\c@thm\makeatother      

\makeatletter\let\c@claim\c@thm\makeatother

\makeatletter\let\c@conjecture\c@thm\makeatother

\newtheorem*{unnumberedtheorem}{Theorem}
\theoremstyle{definition}

\newtheorem{defn}{Definition}[section]
\makeatletter\let\c@defn\c@thm\makeatother
\newtheorem{const}{Construction}[section]
\makeatletter\let\c@const\c@thm\makeatother
\newtheorem{notn}{Notation}[section]
\makeatletter\let\c@notn\c@thm\makeatother

\makeatletter\let\c@convention\c@thm\makeatother

\theoremstyle{remark}

\newtheorem{rmk}{Remark}[section]
\makeatletter\let\c@rmk\c@thm\makeatother

\makeatletter\let\c@ex\c@thm\makeatother

\makeatletter\let\c@observation\c@thm\makeatother

\makeatletter\let\c@warning\c@thm\makeatother

\makeatletter\let\c@digression\c@thm\makeatother

\makeatletter\let\c@answ\c@thm\makeatother

\makeatletter\let\c@answ\c@thm\makeatother

\makeatletter\let\c@aside\c@thm\makeatother

\makeatletter
\let\c@equation\c@thm
\numberwithin{equation}{section}
\makeatother

\crefname{lem}{Lemma}{Lemmas}
\crefname{thm}{Theorem}{Theorems}
\crefname{defn}{Definition}{Definitions}
\crefname{notn}{Notation}{Notations}
\crefname{const}{Construction}{Constructions}
\crefname{prop}{Proposition}{Propositions}
\crefname{rmk}{Remark}{Remarks}
\crefname{cor}{Corollary}{Corollaries}
\crefname{equation}{Display}{Displays}
\crefname{ex}{Example}{Examples}
\crefname{thmalph}{Theorem}{Theorems}
\crefname{answ}{Answer}{Answers}

\newcommand{\cA}{\mathcal{A}}

\newcommand{\cC}{\mathcal{C}}
\newcommand{\cD}{\mathcal{D}}

\newcommand{\cK}{\mathcal{K}}
\newcommand{\cL}{\mathcal{L}}
\newcommand{\cM}{\mathcal{M}}
\newcommand{\cN}{\mathcal{N}}

\newcommand{\cS}{\mathcal{S}}

\newcommand{\cV}{\mathcal{V}}

\newcommand{\bT}{\Theta_2}

\newcommand{\cat}{\cC\!\mathit{at}}
\newcommand{\set}{\cS\!\mathit{et}}
\newcommand{\sset}{\mathit{s}\set}

\newcommand{\psh}[1]{\set^{#1^{\op}}}
\newcommand{\spsh}[1]{\sset^{#1^{\op}}}
\newcommand{\msset}{m\sset}
 \newcommand{\twocat}{2\cat}

 \newcommand{\vcat}[1]{\cat_{#1}}

\DeclareMathOperator{\Ho}{Ho}

\DeclareMathOperator{\Hom}{Hom}
\DeclareMathOperator{\Map}{Map} 

\DeclareMathOperator{\Ob}{Ob}

\DeclareMathOperator{\colim}{colim}
\DeclareMathOperator{\id}{id}

\newcommand{\aamalg}[1]{\underset{#1}{{\amalg}}} %vo: for pushouts
\DeclareMathOperator{\op}{op}

\DeclareMathOperator{\heq}{heq}

\newcommand{\eqDelta}{\Delta[3]_{\mathrm{eq}}}

\newcommand{\Ntheta}{N^{\Theta_2}}

\newcommand{\NRS}{N^{\operatorname{RS}}}

\newcommand{\Deltaop}{\Delta^{\op}}

\newcommand\simpfverA[9]{%
    \def\tempa{#1}%
    \def\tempb{#2}%
    \def\tempc{#3}%
    \def\tempd{#4}%
    \def\tempe{#5}%
    \def\tempf{#6}%
    \def\tempg{#7}%
    \def\temph{#8}%
    \def\tempi{#9}%
}
\newcommand\simpfverAcontinued[5]{%
    \begin{tikzcd}[column sep=1.3cm, row sep=0.4cm, baseline=(current  bounding  box.center), ampersand replacement=\&]
 #5 \& #4 \arrow[l, "\tempc" swap] \&[-7mm] \&[-7mm] #5 \& #4 \arrow[l, "\tempc" swap]\\
 {} \& \&[-7mm] {=} \&[-7mm] \& \\
 #2 \arrow[r, "\tempa" swap]\arrow[uu, "\tempd", ""{name=a031,inner sep=2pt, swap} ] \arrow[uur, "\tempe", ""{name=a02,inner sep=2pt, swap}]\& #3 \arrow[uu, "\tempb" swap]\&[-7mm] \&[-7mm] #2 \arrow[uu, "\tempd", ""{name=a032,inner sep=2pt, swap}]\arrow[r, "\tempa" swap]\& #3\arrow[uu, "\tempb" swap] \arrow[uul, "\tempf", ""{name=a13,inner sep=2pt, swap}]
 \arrow[phantom, from=a031, to=1-2, scalearrow={0.63*0.8}{start023}{end023}]
 \arrow[Rightarrow, to path={(start023)--(end023)\tikztonodes}, "\tempg" ]%from=a031, to path={(a031)--([xshift=1cm]\tikztotarget)}
  \arrow[phantom, from=a02, to=3-2, scalearrow={0.8}{start012}{end012}]
 \arrow[Rightarrow, to path={(start012)--(end012)\tikztonodes}, "\temph" swap ]
  \arrow[phantom, from=a032, to=3-5, scalearrow={0.63*0.8}{start013}{end013}]
 \arrow[Rightarrow, to path={(start013)--(end013)\tikztonodes}, "\tempi" swap ]
   \arrow[phantom, from=a13, to=1-5, scalearrow={0.8}{start123}{end123}]
 \arrow[Rightarrow, to path={(start123)--(end123)\tikztonodes}, "#1" ] %this is the arrow length for the other arrows
 \end{tikzcd}
}

\tikzset{%
scalearrow/.style n args={3}{
  decoration={
    markings,
    mark=at position (1-#1)/2*\pgfdecoratedpathlength
      with {\coordinate (#2);},
    mark=at position (1+#1)/2*\pgfdecoratedpathlength
      with {\coordinate (#3);},
    },
  postaction=decorate,
  }
}
  
\newcommand\simpfverArrow[9]{%
    \def\tempa{#1}%
    \def\tempb{#2}%
    \def\tempc{#3}%
    \def\tempd{#4}%
    \def\tempe{#5}%
    \def\tempf{#6}%
    \def\tempg{#7}%
    \def\temph{#8}%
    \def\tempi{#9}%
}
\newcommand\simpfverArrowcontinued[5]{%
    \begin{tikzcd}[column sep=1.3cm, row sep=0.4cm, baseline=(current  bounding  box.center), ampersand replacement=\&]
 #5 \& #4 \arrow[l, "\tempc" swap] \&[-7mm] \&[-7mm] #5 \& #4 \arrow[l, "\tempc" swap]\\
 {} \& \&[-7mm] {=} \&[-7mm] \& \\
 #2 \arrow[r, "\tempa" swap]\arrow[uu, "\tempd", ""{name=a031,inner sep=2pt, swap} ] \arrow[uur, "\tempe"{description}, ""{name=a02,inner sep=2pt, swap}]\& #3 \arrow[uu, "\tempb" swap]\&[-7mm] \&[-7mm] #2 \arrow[uu, "\tempd", ""{name=a032,inner sep=2pt, swap}]\arrow[r, "\tempa" swap]\& #3\arrow[uu, "\tempb" swap] \arrow[uul, "\tempf"{description}, ""{name=a13,inner sep=2pt, swap}]
 \arrow[phantom, from=a031, to=1-2, scalearrow={0.63*0.8}{start023}{end023}]
 \arrow[Rightarrow, to path={(start023)--(end023)\tikztonodes}, "\tempg" ]%from=a031, to path={(a031)--([xshift=1cm]\tikztotarget)}
  \arrow[phantom, from=a02, to=3-2, scalearrow={0.8}{start012}{end012}]
 \arrow[Rightarrow, to path={(start012)--(end012)\tikztonodes}, "\temph" swap ]
  \arrow[phantom, from=a032, to=3-5, scalearrow={0.63*0.8}{start013}{end013}]
 \arrow[Rightarrow, to path={(start013)--(end013)\tikztonodes}, "\tempi" swap ]
   \arrow[phantom, from=a13, to=1-5, scalearrow={0.8}{start123}{end123}]
 \arrow[Rightarrow, to path={(start123)--(end123)\tikztonodes}, "#1" ] %this is the arrow length for the other arrows
 \end{tikzcd}
}